\documentclass[11pt,a4paper]{article}

\usepackage{amsmath,amssymb,amsfonts,amsthm}
\usepackage{mathtools}
\usepackage{bm}

\usepackage{graphicx}
\usepackage{subfig}

\usepackage{booktabs}
\usepackage{multirow}
\usepackage{array}

\usepackage{geometry}
\usepackage[round,authoryear]{natbib}

\usepackage[hidelinks]{hyperref}

\numberwithin{equation}{section}

\newtheorem{theorem}{Theorem}[section]
\newtheorem{lemma}[theorem]{Lemma}
\newtheorem{corollary}[theorem]{Corollary}
\newtheorem{proposition}[theorem]{Proposition}

\theoremstyle{definition}

\newtheorem{example}[theorem]{Example}

\theoremstyle{remark}
\newtheorem{remark}[theorem]{Remark}

\newenvironment{keywords}
{\par\noindent\textbf{Keywords:}\ }
{\par\vspace{1em}}

\newcommand{\tbl}[1]{%
    \caption{#1}
}

\begin{document}

\begin{center}

{\Large\bfseries Bivariate Rényi--Tsallis Divergence}

\vspace{8pt}

Suramya Surendran\textsuperscript{a}
and
Rajesh G.\textsuperscript{a,*}

\vspace{5pt}

\textsuperscript{a}Department of Statistics, Cochin University of Science and Technology,\\
Cochin 682 022, Kerala, India

\vspace{5pt}

\textit{*Corresponding author: Email: }
\href{mailto:rajeshgstat@gmail.com}{rajeshgstat@gmail.com}

\end{center}

\begin{abstract}
In this paper, a bivariate extension of the differential Rényi–Tsallis divergence is introduced for comparing two absolutely continuous bivariate probability distributions. Fundamental properties of the proposed divergence are studied, including non-negativity, pseudo-additivity under component wise independence, and invariance under monotone and nonsingular linear transformations. Relationships with Rényi, Tsallis, $\alpha$, Bhattacharyya, and Hellinger divergence measures are also established. A sensitivity analysis is performed to compare the Rényi, Tsallis and Rényi-Tsallis divergences under the Marshall-Olkin and Gumbel-type bivariate exponential models, taking the magnitude of Kendall's rank correlation coefficient as a common scale for dependence. The proposed measure is further developed for bivariate order statistics, with particular emphasis on the joint distribution of the minimum and maximum order statistics for several lifetime distributions, including exponential, Weibull, Pareto, and Rayleigh models. When the underlying distributions are unknown, a nonparametric estimator is constructed using bivariate Gaussian kernel density estimators with plug-in bandwidth selection. The proposed estimator is applied to a real life dataset
and bootstrap resampling is used to assess estimation uncertainty. The results show the applicability of the proposed divergence for comparing bivariate distributions and indicate that its finite-sample performance depends on the order parameter and the distributions being compared.
\end{abstract}


\section{Introduction}

The notion of entropy was first introduced by \cite{clausius1865} as a measure of disorder in thermodynamics. This concept was later extended to information theory by \cite{shannon1948} to measure the average uncertainty in a random variable. For a continuous non-negative random variable \( X \) with probability density function \( f_X(x) \), the differential entropy is defined as,  
\begin{equation}
\label{eq:diff_entropy}
H(X)
= - \int_{0}^{\infty} f(x)\, \log f(x)\, dx 
= - E_{X} \left[ \log f(X) \right].
\end{equation}
Among the several generalizations of entropy measures, Rényi entropy introduced by \cite{renyi1961} is one of the most important measures. For a continuous random variable $X$ with probability density function $f_X(x)$ the Rényi entropy of order $\alpha$ is given by,
\begin{equation}
\label{eq:RényiEntropy}
H_{\alpha}(X)
= \frac{1}{1-\alpha}
\log \left(
\int_{0}^{\infty}
f(x)^{\alpha} \, dx
\right) 
= \frac{1}{1-\alpha}
\log \left(
E_{X} \left[ f(X)^{\alpha-1} \right]
\right),
\ \alpha > 0,\ \alpha \neq 1.
\end{equation}
\cite{tsallis1988} introduced another generalization called Tsallis entropy of order $\alpha$. Tsallis entropy is widely used for analyzing non extensive data. For a continuous non-negative random variable $X$ with probability density function $f_X(x)$, the Tsallis entropy is defined as follows,
\begin{equation}
\label{eq:TsallisEntropy}
T_{\alpha}(X)
= \frac{1}{\alpha - 1}
\left[
1 - \int_{0}^{\infty} f(x)^{\alpha} \, dx
\right] 
= \frac{1}{\alpha - 1}
\left[
1 - E_{X}\left( f(X)^{\alpha-1} \right)
\right],
\ \alpha > 0,\ \alpha \neq 1.
\end{equation}
As $\alpha \to 1$, both equation \eqref{eq:RényiEntropy} and \eqref{eq:TsallisEntropy} converge to equation \eqref{eq:diff_entropy}.

Apart from entropy measures \cite{kullback1951} introduced
the Kullback--Leibler (KL) divergence, provides a useful measure for comparing two probability distributions $F_X(\cdot)$ and $G_Y(\cdot)$, and is denoted by $d_{\mathrm{KL}}(F_{\mathbf{X}},G_{\mathbf{Y}})$. It is a measure of the inefficiency of assuming that the distribution is  $g_Y(x)$ when the true distribution is $f_X(x)$. Let $X$ and $Y$ be two continuous non-negative random variables with probability density functions $f_X(x)$ and $g_Y(x)$, respectively. Then the Kullback--Leibler (KL) divergence is defined as,
\begin{equation}
d_{KL}(X,Y)
=
\int_{0}^{\infty}
f(x)\log\frac{f(x)}{g(x)}\,dx 
=
E_{\mathbf{X}}
\left[
\log\left(\frac{f(X)}{g(X)}\right)
\right].
\label{eq:kl}
\end{equation}
\cite{Ebrahimi2007} developed the bivariate case of Kullback--Leibler divergence and studied its properties.
Let $\mathbf{X} = (X_1,X_2)$ and $\mathbf{Y} = (Y_1,Y_2)$ be two  non-negative bivariate random vectors with absolutely continuous joint pdfs $f_{\mathbf{X}}(x_1,x_2)$ and $g_{\mathbf{Y}}(x_1,x_2)$ and the corresponding cdfs $F_{\mathbf{X}}(x_1,x_2)$ and $G_{\mathbf{Y}}(x_1,x_2)$ respectively.  
Then the bivariate KL divergence is defined as,
\begin{equation}
\label{eq:KL_bivariate}
\begin{aligned}
d_{\mathrm{KL}}(F_{\mathbf{X}},G_{\mathbf{Y}})
&= \int_{0}^{\infty}\!\!\int_{0}^{\infty}
f_{\mathbf{X}}(x_1,x_2)
\log\!\left(
\frac{f_{\mathbf{X}}(x_1,x_2)}{g_{\mathbf{Y}}(x_1,x_2)}
\right)
\, dx_1\, dx_2 \\
&= E_{\mathbf{X}}
\left[
\log\!\left(
\frac{f_{\mathbf{X}}(X_1,X_2)}{g_{\mathbf{Y}}(X_1,X_2)}
\right)
\right].
\end{aligned}
\end{equation}
For more information one can refer \cite{contrerasreyes2012}, \cite{zhang2023properties}, \cite{rafflesia2026}.
Several generalizations of the KL divergence exist in the literature. Among them, the Rényi and Tsallis divergence measures of order $\alpha$ are the most studied. The Rényi divergence of order $\alpha$ between two continuous bivariate distributions $F_{\mathbf{X}}$ and $G_{\mathbf{Y}}$ is defined as
\begin{equation}
\label{eq:Rényi}
\begin{aligned}
d^{\alpha}_{R}(F_{\mathbf{X}},G_{\mathbf{Y}})
&=
\frac{1}{\alpha-1}
\left[
\log
\int_{0}^{\infty}\!\!\int_{0}^{\infty}
f_{\mathbf{X}}^{\alpha}(x_1,x_2)
g_{\mathbf{Y}}^{1-\alpha}(x_1,x_2)
\, dx_1\, dx_2 
\right]\\
&=
\frac{1}{\alpha-1}
\left[
\log
E_{\mathbf{X}}
\left(
\frac{f_{\mathbf{X}}(X_1,X_2)}{g_{\mathbf{Y}}(X_1,X_2)}
\right)^{\alpha-1}
\right],
\qquad 0 < \alpha \neq 1 < \infty .
\end{aligned}
\end{equation}
Similarly, the Tsallis divergence of order $\alpha$ is defined as
\begin{equation}
\label{eq:Tsallis}
\begin{aligned}
d^{\alpha}_{T}(F_{\mathbf{X}},G_{\mathbf{Y}})
&=
\frac{1}{\alpha-1}
\left[
\int_{0}^{\infty}\!\!\int_{0}^{\infty}
f_{\mathbf{X}}^{\alpha}(x_1,x_2)
g_{\mathbf{Y}}^{1-\alpha}(x_1,x_2)
\, dx_1\, dx_2
- 1
\right] \\
&=
\frac{1}{\alpha-1}
\left[
E_{\mathbf{X}}
\left(
\frac{f_{\mathbf{X}}(X_1,X_2)}{g_{\mathbf{Y}}(X_1,X_2)}
\right)^{\alpha-1}
- 1
\right],
\ 0 < \alpha \neq 1 < \infty .
\end{aligned}
\end{equation}
It can be shown that, as $\alpha \to 1$, both the Rényi and Tsallis divergence measures converge to the Kullback–Leibler divergence given in \eqref{eq:KL_bivariate}. Furthermore, when $\alpha = 0$, both divergence measures reduce to zero, that is,
$d^{\alpha}_{R}(F_{\mathbf{X}},G_{\mathbf{Y}})=d^{\alpha}_{T}(F_{\mathbf{X}},G_{\mathbf{Y}})=0$.
By combining the concepts of Rényi entropy and Tsallis entropy, \cite{wondie2017joint}  introduced a new measure called the Rényi--Tsallis entropy of order $\alpha$. This generalized entropy includes the Shannon entropy as a special case when the parameter $\alpha$ approaches 1.
Inspired by this, \cite{Khammar2022} introduced the differential Rényi–Tsallis divergence. This new measure combines ideas from Rényi and Tsallis divergences. Let $X$ and $Y$ be two continuous random variables with distribution
functions $F_X(x)$ and $G_Y(x)$, and corresponding probability density
functions $f_X(x)$ and $g_Y(x)$, respectively. The Rényi--Tsallis
divergence between $F_X$ and $G_Y$ is defined as,

\begin{equation}
\label{eq:URT}
d^{\alpha}_{R-T}(F_{\mathbf{X}},G_{\mathbf{Y}})
=
\begin{cases}
\displaystyle
\frac{1}{\alpha - \alpha^{-1}}
\Bigg[
\log
\int_{0}^{\infty}\!
f_{\mathbf{X}}(x)^{\alpha}
g_{\mathbf{Y}}(x)^{1-\alpha}
\, dx \\[2mm]
\hspace{0.8em}
+
\int_{0}^{\infty}\!
f_{\mathbf{X}}^{\alpha}(x)
g_{\mathbf{Y}}^{1-\alpha}(x)
\, dx
- 1
\Bigg],
& 0 < \alpha \neq 1 < \infty, \\[2mm]

\displaystyle
\int_{0}^{\infty}\!
f_{\mathbf{X}}(x)
\log\!\left(
\frac{f_{\mathbf{X}}(x)}{g_{\mathbf{Y}}(x)}
\right)
\, dx,
& \alpha = 1 .
\end{cases}
\end{equation}

\begin{equation*}
\label{eq:URTexp}
d^{\alpha}_{R-T}(F_{\mathbf{X}},G_{\mathbf{Y}})
=
\left\{
\begin{alignedat}{2}
&\frac{1}{\alpha - \alpha^{-1}}
\Bigg[
\log\left[
E_{\mathbf{X}}
\left(
\frac{g_{\mathbf{Y}}(X)}{f_{\mathbf{X}}(X)}
\right)^{1-\alpha}\right]
\\
&\qquad +
E_{\mathbf{X}}
\left(
\frac{g_{\mathbf{Y}}(X)}{f_{\mathbf{X}}(X)}
\right)^{1-\alpha}
- 1
\Bigg],
&\quad 0<\alpha\neq1<\infty,
\\[1mm]
&E_{\mathbf{X}}
\left[
\log\!\left(
\frac{f_{\mathbf{X}}(X)}{g_{\mathbf{Y}}(X)}
\right)
\right],
&\ \alpha=1.
\end{alignedat}
\right.
\end{equation*}
Its major properties were studied and also applied to exponential family distributions and order statistics.

Although Rényi and Tsallis divergences have been extensively studied in multivariate settings, the differential Rényi–Tsallis divergence has received comparatively limited attention beyond the univariate framework. In particular, its bivariate properties, behavior for bivariate order statistics, and nonparametric estimation has not been developed in a unified framework. This motivates the present study.

Rest of this article is organized as follows. In \hyperref[sec:bivariate-rt]{Section~\ref*{sec:bivariate-rt}} we introduce a bivariate differential Rényi–Tsallis divergence and discuss its fundamental properties, including non-negativity, pseudo-additivity under independence, and invariance under monotone and nonsingular linear transformations. We also obtain some bounds for the bivariate Rényi--Tsallis divergence measure. In \hyperref[sec:order-statistics]{Section~\ref*{sec:order-statistics}} we develop Rényi--Tsallis divergence for bivariate order statistics. In \hyperref[sec:nonparametric]{Section~\ref*{sec:nonparametric}}, we study the non parametric estimation of the proposed measure and apply it to a real life dataset. The major conclusions are discussed in \hyperref[sec:conclusion]{Section~\ref*{sec:conclusion}}.

\section{Bivariate Rényi--Tsallis Divergence Measure}
\label{sec:bivariate-rt}

Let $\mathbf{X}=(\mathbf{X_1},\mathbf{X_2})$ and $\mathbf{Y}=(\mathbf{Y_1},\mathbf{Y_2})$ be two absolutely continuous non-negative random vectors having cumulative distribution functions $F_{\mathbf{X}}(x_1,x_2)$ and $G_{\mathbf{Y}}(x_1,x_2)$ and joint probability
density functions $f_{\mathbf{X}}(x_1,x_2)$ and $g_{\mathbf{Y}}(x_1,x_2)$,
respectively. For an order parameter $\alpha>0$, $\alpha\neq1$, define

\[I_\alpha=\int_0^\infty\int_0^\infty
f_X^\alpha(x_1,x_2)
g_Y^{1-\alpha}(x_1,x_2)
,dx_1,dx_2.
\]
Whenever the integral is well-defined, the bivariate differential Rényi–Tsallis divergence of order $\alpha$ is given by,
\begin{equation}
\label{eq:RT}
d^{\alpha}_{R-T}(F_{\mathbf{X}},G_{\mathbf{Y}})
=
\frac{1}{\alpha - \alpha^{-1}}
\Bigg[
\log
I_\alpha
+
I_\alpha
- 1
\Bigg],
\  \alpha >0, \alpha \neq 1 , \\[1mm]
\end{equation}
At $\alpha=1$, the divergence is defined as the Kullback–Leibler divergence:
\begin{equation}
d^{\alpha}_{R-T}(F_{\mathbf{X}},G_{\mathbf{Y}})
=
\int_{0}^{\infty}\!\!\int_{0}^{\infty}
f_{\mathbf{X}}(x_1,x_2)
\log\!\left(
\frac{f_{\mathbf{X}}(x_1,x_2)}{g_{\mathbf{Y}}(x_1,x_2)}
\right)
\, dx_1\, dx_2,
\quad \alpha = 1 .
\end{equation}
provided the suitable support and integrability conditions are satisfied. The same divergence can be written in expectation form. 
\begin{equation}
\label{eq:RTexp}
d^{\alpha}_{R-T}(F_{\mathbf{X}},G_{\mathbf{Y}})
=
\left\{
\begin{alignedat}{2}
&\frac{1}{\alpha - \alpha^{-1}}
\Bigg[
\log
E_{\mathbf{X}}
\left(
\frac{g_{\mathbf{Y}}(X_1,X_2)}{f_{\mathbf{X}}(X_1,X_2)}
\right)^{1-\alpha}
\\
&\qquad +
E_{\mathbf{X}}
\left(
\frac{g_{\mathbf{Y}}(X_1,X_2)}{f_{\mathbf{X}}(X_1,X_2)}
\right)^{1-\alpha}
- 1
\Bigg],
\quad 0<\alpha\neq1<\infty,
\\[2mm]
&E_{\mathbf{X}}
\left[
\log\!\left(
\frac{f_{\mathbf{X}}(X_1,X_2)}{g_{\mathbf{Y}}(X_1,X_2)}
\right)
\right],
\quad \alpha=1.
\end{alignedat}
\right.
\end{equation}
where $E_{\mathbf{X}}$ denotes the expectation with respect to $f_X$. The proposed formulation is a mixture of the logarithmic and linear components of the power-integral functional. This establishes its connection between the Renyi and Tsallis divergences. Using \eqref{eq:Rényi}, \eqref{eq:Tsallis}, and \eqref{eq:RT},
the differential Rényi--Tsallis divergence can be represented as
\begin{equation}
\label{eq:RTrelation}
d^{\alpha}_{R\!-\!T}(F_{\mathbf{X}},G_{\mathbf{Y}})
=
\frac{\alpha}{\alpha+1}
\left[
d^{\alpha}_{R}(F_{\mathbf{X}},G_{\mathbf{Y}})
+
d^{\alpha}_{T}(F_{\mathbf{X}},G_{\mathbf{Y}})
\right].
\end{equation}
As $\alpha \to 1$, equation \eqref{eq:RTrelation} reduces to the measure given in \eqref{eq:KL_bivariate}. Moreover, \eqref{eq:RTrelation} indicates that the behavior of the divergence measure $d_{R-T}^{\alpha}(F_{\mathbf{X}}, G_{\mathbf{Y}})$ is similar to that of $d_R^\alpha(F_{\mathbf{X}}, G_{\mathbf{Y}})$ and $d_T^\alpha(F_{\mathbf{X}}, G_{\mathbf{Y}})$ with respect to changes in $\alpha$. 
\begin{remark}
    Let  $I_\alpha(F_X,G_Y) = \int_{0}^{\infty}\!\!\int_{0}^{\infty}
f_{\mathbf{X}}^{\alpha}(x_1,x_2)
g_{\mathbf{Y}}^{1-\alpha}(x_1,x_2)
\, dx_1\, dx_2$. Now by interchanging $F_X$ and $G_Y$ we get $I_\alpha(G_Y,F_X) = \int_{0}^{\infty}\!\!\int_{0}^{\infty}
g_{\mathbf{Y}}^{\alpha}(x_1,x_2)
f_{\mathbf{X}}^{1-\alpha}(x_1,x_2)
\, dx_1\, dx_2$. In general these two integrals are not equal. Hence, \[d^\alpha_{R-T}(F_X,G_Y)
\neq
d^\alpha_{R-T}(G_Y,F_X).\]
That is, the bivariate Rényi--Tsallis divergence is generally asymmetric. When
$\alpha=\tfrac{1}{2}$, we have
$
I_{\frac{1}{2}}(F_X,G_Y)
=
\int_{0}^{\infty}\!\!\int_{0}^{\infty}
f_{\mathbf{X}}^{\frac{1}{2}}(x_1,x_2)
g_{\mathbf{Y}}^{\frac{1}{2}}(x_1,x_2)
\,dx_1\,dx_2 
=
I_{\frac{1}{2}}(G_Y,F_X).
$
Hence, the bivariate Rényi--Tsallis divergence becomes symmetric at
$\alpha=\tfrac{1}{2}$ since the underlying integral reduces to the Bhattacharyya coefficient, which is invariant under the interchange of the two distributions. 
\end{remark}

\subsection*{Marginal--conditional representation}

Let $X=(X_1,X_2)$ and $Y=(Y_1,Y_2)$ be absolutely continuous
non-negative bivariate random vectors with joint pdfs
$
f_X(x_1,x_2)
=
f_{X_1}(x_1)f_{X_2\mid X_1}(x_2\mid x_1)
$
and
$
g_Y(x_1,x_2)
=
g_{Y_1}(x_1)g_{Y_2\mid Y_1}(x_2\mid x_1).
$
Using these representations, the integral in the definition
of the bivariate R\'enyi--Tsallis divergence can be written as
\begin{equation}
\begin{aligned}
I_\alpha(F_X,G_Y)
&=
\int_0^\infty\int_0^\infty
f_X^\alpha(x_1,x_2)
g_Y^{1-\alpha}(x_1,x_2)
\,dx_1\,dx_2\\
&=
\int_0^\infty
f_{X_1}^{\alpha}(x_1)
g_{Y_1}^{1-\alpha}(x_1)
\mathcal{I}_\alpha(x_1)\,dx_1.
\label{eq:int}
\end{aligned}
\end{equation}
where
\[
\mathcal{I}_\alpha(x_1)
=
\int_0^\infty
f_{X_2\mid X_1}^{\alpha}(x_2\mid x_1)
g_{Y_2\mid Y_1}^{1-\alpha}(x_2\mid x_1)
\,dx_2.
\]
Then, the bivariate R\'enyi--Tsallis divergence can be
expressed as
\begin{equation}
d_{R-T}^{\alpha}(F_X,G_Y)
=
\frac{1}{\alpha-\alpha^{-1}}
\left[
\log I_\alpha(F_X,G_Y)
+
I_\alpha(F_X,G_Y)-1
\right],
\tag{15}
\end{equation}
where $I_\alpha(F_X,G_Y)$ is given by  \eqref{eq:int}.

This representation separates the contribution of the marginal densities of $X_1$ and $Y_1$ from the integral containing the conditional densities of $X_2$ given $X_1$ and $Y_2$ given $Y_1$. It thus provides a convenient way to study how the bivariate divergence depends on marginal and conditional distributional information.  However, it does not imply an additive decomposition of the divergence into
marginal and conditional parts, in general.

\begin{example}
Let \(X=(X_1,X_2)\) and \(Y=(Y_1,Y_2)\) follow bivariate Pareto Type--II distributions with probability density functions
$
f_X(x_1,x_2)
=
\theta_1(\theta_1+1)(1+x_1+x_2)^{-(\theta_1+2)},
$
and
$
g_Y(x_1,x_2)
=
\theta_2(\theta_2+1)(1+x_1+x_2)^{-(\theta_2+2)},
$
where \(x_1,x_2>0\), \(\theta_1>0\), and \(\theta_2>0\).
The Rényi--Tsallis divergence measure between \(F_X\) and \(G_Y\) is given by
\begin{equation}
d_{R\text{-}T}^{\alpha}(F_X,G_Y)
=
\left\{
\begin{aligned}
&\frac{1}{\alpha-\alpha^{-1}}
[\log I_{\alpha}
+
\left(I_{\alpha}-1\right)],
\qquad
\alpha>0,\ \alpha\neq1,
\\[1em]
&\frac{
(\theta_2-\theta_1)(1+2\theta_1)
}{
\theta_1(1+\theta_1)
}
+
\log
\left(
\frac{
\theta_1(1+\theta_1)
}{
\theta_2(1+\theta_2)
}
\right),
\qquad
\alpha = 1.
\end{aligned}
\right.
\end{equation}
where
\[
I_{\alpha}
=
\frac{
\left[\theta_1(\theta_1+1)\right]^{\alpha}
\left[\theta_2(\theta_2+1)\right]^{1-\alpha}
}{
\left(\alpha\theta_1+(1-\alpha)\theta_2+1\right)
\left(\alpha\theta_1+(1-\alpha)\theta_2\right)
}.
\]
For \(\alpha=0.5\), the values of the bivariate Rényi--Tsallis divergence for different values of \(\theta_1\) and \(\theta_2\) are presented in Table~\ref{tab:rt_divergence}.
\end{example}

\begin{table}
\centering
\tbl{Bivariate Rényi--Tsallis divergence values for different values of $\theta_1$ and $\theta_2$}
{\scriptsize
\setlength{\tabcolsep}{3pt}
\renewcommand{\arraystretch}{0.9}
\begin{tabular}{lcccccccccc}
\toprule
$\theta_1 \backslash \theta_2$
& 1 & 2 & 3 & 4 & 5 & 6 & 7 & 8 & 9 & 10 \\
\midrule
1
& 0.000 & 0.104 & 0.257 & 0.401 & 0.528 & 0.640 & 0.738 & 0.826 & 0.905 & 0.977 \\
2
& 0.104 & 0.000 & 0.041 & 0.119 & 0.206 & 0.292 & 0.373 & 0.449 & 0.521 & 0.587 \\
3
& 0.257 & 0.041 & 0.000 & 0.022 & 0.069 & 0.127 & 0.188 & 0.249 & 0.309 & 0.366 \\
4
& 0.401 & 0.119 & 0.022 & 0.000 & 0.014 & 0.046 & 0.087 & 0.132 & 0.179 & 0.226 \\
5
& 0.528 & 0.206 & 0.069 & 0.014 & 0.000 & 0.009 & 0.032 & 0.063 & 0.098 & 0.135 \\
6
& 0.640 & 0.292 & 0.127 & 0.046 & 0.009 & 0.000 & 0.007 & 0.024 & 0.048 & 0.076 \\
7
& 0.738 & 0.373 & 0.188 & 0.087 & 0.032 & 0.007 & 0.000 & 0.005 & 0.019 & 0.038 \\
8
& 0.826 & 0.449 & 0.249 & 0.132 & 0.063 & 0.024 & 0.005 & 0.000 & 0.004 & 0.015 \\
9
& 0.905 & 0.521 & 0.309 & 0.179 & 0.098 & 0.048 & 0.019 & 0.004 & 0.000 & 0.003 \\
10
& 0.977 & 0.587 & 0.366 & 0.226 & 0.135 & 0.076 & 0.038 & 0.015 & 0.003 & 0.000 \\
\bottomrule
\end{tabular}}
\label{tab:rt_divergence}
\end{table}

\begin{remark}
Let \(X=(X_1,X_2)\), \(Y=(Y_1,Y_2)\), and \(Z=(Z_1,Z_2)\) be Pareto Type--II random vectors with parameters
$
\theta_1=1,\ \theta_2=5,\ \theta_3=10
$ respectively.
To verify the triangular inequality property of the Rényi--Tsallis divergence, consider
\[
d^{{\tfrac{1}{2}}}_{R-T}(X,Y)
=
0.528,\ d^{{\tfrac{1}{2}}}_{R-T}(Y,Z)
=
0.135,\ d^{{\tfrac{1}{2}}}_{R-T}(X,Z)
=
0.976.
\]
\[
d^{{\tfrac{1}{2}}}_{R-T}(X,Y)+d^{\tfrac{1}{2}}_{R-T}(Y,Z)
=
0.528158+0.135348
=
0.663506.
\]
\[
d^{{\tfrac{1}{2}}}_{R-T}(X,Z)
=
0.976563
>
0.663506
=
d^{{\tfrac{1}{2}}}_{R-T}(X,Y)+d^{\tfrac{1}{2}}_{R-T}(Y,Z),
\]
\[
d^{\alpha}_{R-T}(X,Z)
\nleq
d^{\alpha}_{R-T}(X,Y)+d^{\alpha}_{R-T}(Y,Z).
\]
Hence, the triangular inequality is not satisfied. 
\end{remark}

\begin{remark}
As discussed earlier, in general the Rényi--Tsallis divergence is not symmetric. However, for $\alpha=\tfrac{1}{2}$, the divergence becomes symmetric. This property is clearly visible in Table~\ref{tab:rt_divergence}, in which the divergence value for each pair $(\theta_1,\theta_2)$ is equal to that of $(\theta_2,\theta_1)$. For example,
\[
d_{R-T}^{{\tfrac{1}{2}}}(1,4)=d_{R-T}^{{\tfrac{1}{2}}}(4,1)=0.401
\]
and
\[
d_{R-T}^{{\tfrac{1}{2}}}(2,7)=d_{R-T}^{{\tfrac{1}{2}}}(7,2)=0.373.
\]
Thus, the symmetry in Table~\ref{tab:rt_divergence} is due to the choice $\alpha=\tfrac{1}{2}$, but the Rényi--Tsallis divergence is not symmetric in general. So, from remark 1 and 2 it is not a metric. The following theorem proves that the bivariate Rényi--Tsallis divergence is non--negative.
\end{remark}

\subsection*{Non--negativity}
\begin{theorem}
Let $\mathbf{X}=(X_1,X_2)$ and $\mathbf{Y}=(Y_1,Y_2)$ be absolutely continuous bivariate random vectors with cdfs $F_X(x_1,x_2)$ and $G_Y(x_1,x_2)$, and corresponding joint pdfs $f_X(x_1,x_2)$ and $g_Y(x_1,x_2)$, respectively. Then
$
d^\alpha_{R-T}(F_X,G_Y)\geq 0,
\ for \ $\mbox{$0<\alpha\neq1<\infty$}$.
$

\end{theorem}

\begin{proof}

Let
\[
I_\alpha(F_X,G_Y)=
\int_0^\infty\int_0^\infty
f_X^\alpha(x_1,x_2)
g_Y^{1-\alpha}(x_1,x_2)
\,dx_1\,dx_2.
\]
The Rényi--Tsallis divergence is defined as
\[
d^\alpha_{R-T}(F_X,G_Y)
=
\frac{1}{\alpha-\alpha^{-1}}
\left[
\log I_\alpha+I_\alpha-1
\right].
\]
Since $f_X(x_1,x_2)$ is a probability density, we can write
$
I_\alpha(F_X,G_Y)
=
\int_0^\infty\int_0^\infty
f_X(x_1,x_2)
\left[
\frac{g_Y(x_1,x_2)}
{f_X(x_1,x_2)}
\right]^{1-\alpha}
\,dx_1\,dx_2.
$
Define
$
Z=
\frac{g_Y(X_1,X_2)}
{f_X(X_1,X_2)}.
$
Then
$
I_\alpha(F_X,G_Y)
=
E_X\left[Z^{1-\alpha}\right].
$
Moreover,
\[
\begin{aligned}
E_X(Z)
&=
\int_0^\infty\int_0^\infty
f_X(x_1,x_2)
\frac{g_Y(x_1,x_2)}
{f_X(x_1,x_2)}
\,dx_1\,dx_2\\
&=
\int_0^\infty\int_0^\infty
g_Y(x_1,x_2)
\,dx_1\,dx_2=1.
\end{aligned}
\]
We consider two cases.
\\
\textbf{Case 1: $0<\alpha<1$.}
In this case,
$
0<1-\alpha<1.
$
Hence, the function
$
\phi(z)=z^{1-\alpha},\ z>0,
$
is concave. By Jensen's inequality for multivariate case (\cite{perlman1974jensen}),
\[
I_\alpha(F_X,G_Y)
=
E_X[\phi(Z)]
\leq
\phi(E_X[Z]).
\]
Since $E_X(Z)=1$, we obtain
$
I_\alpha(F_X,G_Y)\leq\phi(1)=1.
$
Therefore,
$
\log I_\alpha(F_X,G_Y)\leq0
$
and
$
I_\alpha(F_X,G_Y)-1\leq0.
$
Thus,
$
\log I_\alpha(F_X,G_Y)+I_\alpha-1\leq0.
$
For $0<\alpha<1$,
$
\alpha-\alpha^{-1}<0,
$
and hence
$
\frac{1}{\alpha-\alpha^{-1}}<0.
$
Therefore,
\[
d^\alpha_{R-T}(F_X,G_Y)
=
\frac{1}{\alpha-\alpha^{-1}}
\left[
\log I_\alpha+I_\alpha-1
\right]
\geq0.
\\
\]
\textbf{Case 2: $\alpha>1$.}
In this case,
$
1-\alpha<0.
$
Thus,
$
\phi(z)=z^{1-\alpha}
$
is convex on $(0,\infty)$. By Jensen's inequality for multivariate case,
\[
I_\alpha(F_X,G_Y)
=
E_X[\phi(Z)]
\geq
\phi(E_X[Z]).
\]
Since $E_X(Z)=1$,
$
I_\alpha(F_X,G_Y)\geq\phi(1)=1.
$
Consequently,
$
\log I_\alpha(F_X,G_Y)\geq0
$
and
$
I_\alpha(F_X,G_Y)-1\geq0.
$
Hence,
$
\log I_\alpha(F_X,G_Y)+I_\alpha(F_X,G_Y)-1\geq0.
$
Also, for $\alpha>1$,
$
\alpha-\alpha^{-1}>0,
$
so that
$
\frac{1}{\alpha-\alpha^{-1}}>0.
$
Therefore,
\[
d^\alpha_{R-T}(F_X,G_Y)\geq0.
\]
The equality condition holds when $F_{\mathbf{X}}=G_{\mathbf{Y}}$, that is, when the two bivariate random vectors have the same distribution function almost everywhere.
This completes the proof.

\end{proof}
The behavior of these divergence measures as $\alpha$ varies is given in the following lemma.
\begin{lemma}
\label{lem:divergence_behavior}
The divergence measures $d_R^\alpha(F_{\mathbf{X}},G_{\mathbf{Y}})$, $d_T^\alpha(F_{\mathbf{X}},G_{\mathbf{Y}})$, and $d_{R-T}^\alpha(F_{\mathbf{X}},G_{\mathbf{Y}})$ exhibit the following behavior based on $\alpha$,

\begin{itemize}
\item If $\alpha = 0$, then 
$d_R^{\alpha}(F_{\mathbf{X}},G_{\mathbf{Y}})=0$, $d_T^{\alpha}(F_{\mathbf{X}},G_{\mathbf{Y}})=0$, and $d_{R-T}^{\alpha}(F_{\mathbf{X}},G_{\mathbf{Y}})=0$.

\item If $0<\alpha<1$, then 
$d_R^{\alpha}(F_{\mathbf{X}},G_{\mathbf{Y}})\ge 0$, $d_T^{\alpha}(F_{\mathbf{X}},G_{\mathbf{Y}})\ge 0$, and $d_{R-T}^{\alpha}(F_{\mathbf{X}},G_{\mathbf{Y}})\ge 0$.

\item If $\alpha = 1$, then 
$d_R^{\alpha}(F_{\mathbf{X}},G_{\mathbf{Y}})=d_T^{\alpha}(F_{\mathbf{X}},G_{\mathbf{Y}})=d_{R-T}^{\alpha}(F_{\mathbf{X}},G_{\mathbf{Y}})=\mathrm{KL}(F_{\mathbf{X}},G_{\mathbf{Y}})\ge 0$.

\item If $\alpha > 1$, then 
$d_R^{\alpha}(F_{\mathbf{X}},G_{\mathbf{Y}})\ge 0$, $d_T^{\alpha}(F_{\mathbf{X}},G_{\mathbf{Y}})\ge 0$, and $d_{R-T}^{\alpha}(F_{\mathbf{X}},G_{\mathbf{Y}})\ge 0$.
\end{itemize}

\end{lemma}
\subsection*{Pseudo-additivity}
\begin{theorem}
Let $\mathbf{X}=(X_1,X_2)$ and $\mathbf{Y}=(Y_1,Y_2)$ be absolutely continuous non-negative bivariate random vectors with joint pdfs $f_X(x_1,x_2)$ and $g_Y(x_1,x_2)$, respectively. Assume that the components are independent under both distributions, that is,
$
f_X(x_1,x_2)=f_{X_1}(x_1)f_{X_2}(x_2)
$
and
$
g_Y(x_1,x_2)=g_{Y_1}(x_1)g_{Y_2}(x_2).
$
For $\alpha>0$, $\alpha\neq1$, let
$
I_{\alpha,i}
=
\int_0^\infty
f_{X_i}^{\alpha}(x_i)
g_{Y_i}^{1-\alpha}(x_i)
\,dx_i,
\ i=1,2.
$
Then the Rényi--Tsallis divergence satisfies the pseudo-additivity relation
\[
d^\alpha_{R-T}(F_X,G_Y)
=
d^\alpha_{R-T}(F_{X_1},G_{Y_1})
+
d^\alpha_{R-T}(F_{X_2},G_{Y_2})
+
\frac{(I_{\alpha,1}-1)(I_{\alpha,2}-1)}
{\alpha-\alpha^{-1}}.
\]
\end{theorem}

\begin{proof}

Consider
$
I_\alpha(F_X,G_Y)
=
\int_0^\infty\int_0^\infty
f_X^\alpha(x_1,x_2)
g_Y^{1-\alpha}(x_1,x_2)
\,dx_1\,dx_2.
$
As the components are independent under both the distributions,
\[
\begin{aligned}
I_\alpha(F_X,G_Y)
&=
\int_0^\infty\int_0^\infty
f_{X_1}^{\alpha}(x_1)
f_{X_2}^{\alpha}(x_2)
g_{Y_1}^{1-\alpha}(x_1)
g_{Y_2}^{1-\alpha}(x_2)
\,dx_1\,dx_2\\
&=
\left[
\int_0^\infty
f_{X_1}^{\alpha}(x_1)
g_{Y_1}^{1-\alpha}(x_1)
\,dx_1
\right]
\left[
\int_0^\infty
f_{X_2}^{\alpha}(x_2)
g_{Y_2}^{1-\alpha}(x_2)
\,dx_2
\right].
\end{aligned}
\]
Thus,
$
I_\alpha(F_X,G_Y)=I_{\alpha,1}I_{\alpha,2}.
$
By the definition of the Rényi--Tsallis divergence,
\[
d^\alpha_{R-T}(F_X,G_Y)
=
\frac{1}{\alpha-\alpha^{-1}}
\left[
\log I_\alpha(F_X,G_Y)
+
I_\alpha(F_X,G_Y)-1
\right].
\]
Substituting $I_\alpha(F_X,G_Y)=I_{\alpha,1}I_{\alpha,2}$ gives
\[
d^\alpha_{R-T}(F_X,G_Y)
=
\frac{1}{\alpha-\alpha^{-1}}
\left[
\log(I_{\alpha,1}I_{\alpha,2})
+
I_{\alpha,1}I_{\alpha,2}-1
\right].
\]
Using
$
\log(I_{\alpha,1}I_{\alpha,2})
=
\log I_{\alpha,1}+\log I_{\alpha,2}
$
and
$
I_{\alpha,1}I_{\alpha,2}-1
=
(I_{\alpha,1}-1)
+
(I_{\alpha,2}-1)
+
(I_{\alpha,1}-1)(I_{\alpha,2}-1),
$
we obtain
\[
\begin{aligned}
d^\alpha_{R-T}(F_X,G_Y)
&=
\frac{1}{\alpha-\alpha^{-1}}
\left[
\log I_{\alpha,1}+I_{\alpha,1}-1
\right]+
\frac{1}{\alpha-\alpha^{-1}}
\left[
\log I_{\alpha,2}+I_{\alpha,2}-1
\right]\\
&\quad+
\frac{(I_{\alpha,1}-1)(I_{\alpha,2}-1)}
{\alpha-\alpha^{-1}}.
\end{aligned}
\]
Hence,
\[
d^\alpha_{R-T}(F_X,G_Y)
=
d^\alpha_{R-T}(F_{X_1},G_{Y_1})
+
d^\alpha_{R-T}(F_{X_2},G_{Y_2})
+
\frac{(I_{\alpha,1}-1)(I_{\alpha,2}-1)}
{\alpha-\alpha^{-1}}
.
\]
Therefore, the Rényi--Tsallis divergence satisfies a pseudo-additivity property under independence.
\end{proof}

\begin{example}
Let $\mathbf{X}=(X_1,X_2)$ and $\mathbf{Y}=(Y_1,Y_2)$ follow independent bivariate exponential distributions with parameters $(\lambda_1,\lambda_2)$ 
and $(\mu_1,\mu_2)$, respectively. The joint density functions are
$
f_{\mathbf{X}}(x_1,x_2)
=
\lambda_1 \lambda_2 \exp(-\lambda_1 x_1 - \lambda_2 x_2),
g_{\mathbf{Y}}(x_1,x_2)
=
\mu_1 \mu_2 \exp(-\mu_1 x_1 - \mu_2 x_2),
$
for $x_1,x_2>0$.
Then,
\begin{equation*}
I_{\alpha}(F_{\mathbf{X}}, G_{\mathbf{Y}})
=
\int_{0}^{\infty}\!\!\int_{0}^{\infty}
f_{\mathbf{X}}^{\alpha}(x_1,x_2)\,
g_{\mathbf{Y}}^{\,1-\alpha}(x_1,x_2)\,dx_1\,dx_2 
=
\prod_{i=1}^{2}
\frac{\lambda_i^{\alpha}\mu_i^{\,1-\alpha}}
{\alpha \lambda_i + (1-\alpha)\mu_i}.
\end{equation*}
Substituting into \eqref{eq:RT}, we obtain
\begin{equation}
\label{eq:RT_gumbel_bivariate}
d^{\alpha}_{R\!-\!T}(F_{\mathbf{X}}, G_{\mathbf{Y}})
=
\frac{1}{\alpha-\alpha^{-1}}
\left\{
\log
\prod_{i=1}^{2}
\left(
\frac{\lambda_i^{\alpha}\mu_i^{\,1-\alpha}}
{\alpha \lambda_i + (1-\alpha)\mu_i}
\right)
+
\prod_{i=1}^{2}
\frac{\lambda_i^{\alpha}\mu_i^{\,1-\alpha}}
{\alpha \lambda_i + (1-\alpha)\mu_i}
-1
\right\}.
\end{equation}
Further, the marginal divergences are given by
\begin{equation*}
\label{eq:RT_gumbel_marg1}
d^{\alpha}_{R\!-\!T}(F_{X_1}, G_{Y_1})
=
\frac{1}{\alpha-\alpha^{-1}}
\left[
\log
\left(
\frac{\lambda_1^{\alpha}\mu_1^{\,1-\alpha}}
{\alpha \lambda_1 + (1-\alpha)\mu_1}
\right)
+
\frac{\lambda_1^{\alpha}\mu_1^{\,1-\alpha}}
{\alpha \lambda_1 + (1-\alpha)\mu_1}
-1
\right],
\end{equation*}

\begin{equation*}
\label{eq:RT_gumbel_marg2}
d^{\alpha}_{R\!-\!T}(F_{X_2}, G_{Y_2})
=
\frac{1}{\alpha-\alpha^{-1}}
\left[
\log
\left(
\frac{\lambda_2^{\alpha}\mu_2^{\,1-\alpha}}
{\alpha \lambda_2 + (1-\alpha)\mu_2}
\right)
+
\frac{\lambda_2^{\alpha}\mu_2^{\,1-\alpha}}
{\alpha \lambda_2 + (1-\alpha)\mu_2}
-1
\right].
\end{equation*}
Hence,
\begin{equation}
\label{eq:RT_pseudo_deviation}
\begin{aligned}
&d^{\alpha}_{R\!-\!T}(F_{\mathbf{X}}, G_{\mathbf{Y}})
-
\left[
d^{\alpha}_{R\!-\!T}(F_{X_1}, G_{Y_1})
+
d^{\alpha}_{R\!-\!T}(F_{X_2}, G_{Y_2})
\right] \\
&=
\frac{1}{\alpha-\alpha^{-1}}
\left[
\prod_{i=1}^{2}
\frac{\lambda_i^{\alpha}\mu_i^{\,1-\alpha}}
{\alpha \lambda_i + (1-\alpha)\mu_i}
-
\frac{\lambda_1^{\alpha}\mu_1^{\,1-\alpha}}
{\alpha \lambda_1 + (1-\alpha)\mu_1}
-
\frac{\lambda_2^{\alpha}\mu_2^{\,1-\alpha}}
{\alpha \lambda_2 + (1-\alpha)\mu_2}
+1
\right].
\end{aligned}
\end{equation}

The expression \eqref{eq:RT_pseudo_deviation} verifies the pseudo-additive structure of the 
Rényi--Tsallis divergence for the independent bivariate exponential distribution.
\end{example}

\begin{corollary}
Under the assumptions of the pseudo-additivity theorem, suppose that the component integrals $I_{\alpha,1}$ and $I_{\alpha,2}$ are finite and well-defined. Then
\[
d_{R-T}^{\alpha}(F_X,G_Y)
\begin{cases}
\displaystyle
\leq d_{R-T}^{\alpha}(F_{X_1},G_{Y_1})
+d_{R-T}^{\alpha}(F_{X_2},G_{Y_2}),
&0<\alpha<1,\\[6pt]
\displaystyle
\geq d_{R-T}^{\alpha}(F_{X_1},G_{Y_1})
+d_{R-T}^{\alpha}(F_{X_2},G_{Y_2}),
&\alpha>1.
\end{cases}
\]
\end{corollary}
\begin{proof}
By the pseudo-additivity theorem,
\begin{equation}
\label{eq:psead}
d_{R-T}^{\alpha}(F_X,G_Y)
=d_{R-T}^{\alpha}(F_{X_1},G_{Y_1})
+d_{R-T}^{\alpha}(F_{X_2},G_{Y_2})+C_\alpha,
\end{equation}
where
$
C_\alpha=
\frac{(I_{\alpha,1}-1)(I_{\alpha,2}-1)}
{\alpha-\alpha^{-1}}.
$
For \(0<\alpha<1\), the inequalities \(I_{\alpha,i}\leq1\), \(i=1,2\), imply
$
(I_{\alpha,1}-1)(I_{\alpha,2}-1)\geq0.
$
Since \(\alpha-\alpha^{-1}<0\), it follows that \(C_\alpha\leq0\).
For \(\alpha>1\), we have \(I_{\alpha,i}\geq1\), so
$
(I_{\alpha,1}-1)(I_{\alpha,2}-1)\geq0.
$
Since \(\alpha-\alpha^{-1}>0\), we obtain \(C_\alpha\geq0\). Substituting into \eqref{eq:psead} gives the result.
\end{proof}

\subsection*{Invariant under Monotone Transformation}
\begin{theorem}
Let $X =(X_1,X_2)$ and $Y =(Y_1,Y_2)$ be two bivariate random vectors admitting absolutely continuous distribution functions with joint pdfs $f_X(x_1,x_2)$ and $g_Y(x_1,x_2)$, respectively. Consider $S =(S_1,S_2)$ and $T =(T_1,T_2)$ such that $S_i=\phi_i(X_i)$ and $T_i=\phi_i(Y_i)$ for $i=1,2$, where $\phi_i(\cdot)$ are monotone and differentiable. Then
\[
d^{\alpha}_{R\text{-}T}(F_S^*,G_T^*) = d^{\alpha}_{R\text{-}T}(F_X,G_Y).
\]
\end{theorem}

\begin{proof}
\[
\begin{aligned}
d^{\alpha}_{R\text{-}T}(F_X,G_Y)
&= \frac{1}{\alpha-\alpha^{-1}} \Bigg[
\log \iint f_X^\alpha(x_1,x_2)\, g_Y^{\,1-\alpha}(x_1,x_2)\,dx_1\,dx_2 \\
&\quad + \iint f_X^\alpha(x_1,x_2)\, g_Y^{\,1-\alpha}(x_1,x_2)\,dx_1\,dx_2 - 1
\Bigg].
\end{aligned}
\]
Let $f_S^*(s_1,s_2)$ and $g_T^*(s_1,s_2)$ denote the joint pdfs of $S$ and $T$, respectively. Then
$
f_S^*(s_1,s_2) = f_X(x_1,x_2)\,|J|^{-1}, \
g_T^*(s_1,s_2) = g_Y(x_1,x_2)\,|J|^{-1},
$
where $J$ is the Jacobian of the transformation given by
$
J = \left|\frac{\partial(s_1,s_2)}{\partial(x_1,x_2)}\right|
=
\begin{vmatrix}
\frac{\partial s_1}{\partial x_1} & \frac{\partial s_1}{\partial x_2} \\
\frac{\partial s_2}{\partial x_1} & \frac{\partial s_2}{\partial x_2}
\end{vmatrix}.
$ Also $ds_1\,ds_2 = |J|\,dx_1\,dx_2
$
\[
\begin{aligned}
d^{\alpha}_{R\text{-}T}(F_S^*,G_T^*)
&= \frac{1}{\alpha-\alpha^{-1}}\Bigg[\log \iint f_S^{*\alpha}(s_1,s_2)\, g_T^{*1-\alpha}(s_1,s_2)\,ds_1\,ds_2 \\
&\quad + \iint f_S^{*\alpha}(s_1,s_2)\, g_T^{*1-\alpha}(s_1,s_2)\,ds_1\,ds_2 -1 \Bigg] \\
&= \frac{1}{\alpha-\alpha^{-1}}
\Bigg[
\log \iint \big[f_X(x_1,x_2)\,|J|^{-1}\big]^{\alpha}
\big[g_Y(x_1,x_2)\,|J|^{-1}\big]^{1-\alpha}\, ds_1\, ds_2 \\
&\quad + \iint \big[f_X(x_1,x_2)\,|J|^{-1}\big]^{\alpha}
\big[g_Y(x_1,x_2)\,|J|^{-1}\big]^{1-\alpha}\, ds_1\, ds_2 - 1
\Bigg]\\
&= \frac{1}{\alpha-\alpha^{-1}}\Bigg[\log \iint f_X^\alpha(x_1,x_2)\, g_Y^{\,1-\alpha}(x_1,x_2)\,|J|^{-1} ds_1\,ds_2 \\
&\quad + \iint f_X^\alpha(x_1,x_2)\, g_Y^{\,1-\alpha}(x_1,x_2)\,|J|^{-1} ds_1\,ds_2 -1 \Bigg]\\
&= \frac{1}{\alpha-\alpha^{-1}}\Bigg[\log \iint f_X^\alpha(x_1,x_2)\, g_Y^{\,1-\alpha}(x_1,x_2)\,dx_1\,dx_2 \\
&\quad + \iint f_X^\alpha(x_1,x_2)\, g_Y^{\,1-\alpha}(x_1,x_2)\,dx_1\,dx_2 -1 \Bigg]
= d^\alpha_{R\text{-}T}(F_X,G_Y).
\end{aligned}
\]
\end{proof}

\subsection*{Invariance under Linear Transformation}

\begin{corollary}
Let $X=(X_1,X_2)$ and $Y=(Y_1,Y_2)$ be absolutely continuous, bivariate random vectors with joint pdfs $f_X(x_1,x_2)$ and $g_Y(x_1,x_2)$, respectively.  
Consider $S=(S_1,S_2)$ and $T=(T_1,T_2)$ such that
$
S' = AX' + b, \quad T' = AY' + b,
$
where $A$ is a non-singular $2 \times 2$ matrix and $b$ is a $2 \times 1$ constant vector. Then
\[
d^{\alpha}_{R\text{-}T}(F_S^*, G_T^*) = d^\alpha_{R\text{-}T}(F_X, G_Y).
\]
That is the bivariate Rényi--Tsallis divergence is invariant under linear transformation.
\end{corollary}

\begin{proof}
This result follows directly from the Theorem 2.7.
\end{proof}
\subsection*{Relationship with some bivariate divergence measures}
Let $X = (X_1, X_2)$ and $Y = (Y_1, Y_2)$ be two non-negative absolutely continuous bivariate random vectors with cdfs $F_X(x_1, x_2)$ and $G_Y(y_1, y_2)$ and density functions $f_X(x_1, x_2)$ and $g_Y(x_1, x_2)$, respectively. Then,
\[
d_c^\alpha(F_X, G_Y)
=
\frac{1}{\alpha(1-\alpha)}
\left[
1-
\iint
f_X^\alpha(x_1,x_2)
g_Y^{\,1-\alpha}(x_1,x_2)
\,dx_1\,dx_2
\right]
\]
is known as the bivariate $\alpha$-divergence measure introduced by \cite{Cichocki2009}.
From this, we obtain
\[
\iint
f_X^\alpha(x_1,x_2)
g_Y^{\,1-\alpha}(x_1,x_2)
\,dx_1\,dx_2
=
1-\alpha(1-\alpha)d_c^\alpha(F_X,G_Y).
\]
Substituting this into the equation for
$d_{R-T}^\alpha(F_X,G_Y)$, we get
\[
d_{R-T}^\alpha(F_X,G_Y)
=
\frac{1}{\alpha-\alpha^{-1}}
\left\{
\log\left[
1-\alpha(1-\alpha)
d_c^\alpha(F_X,G_Y)
\right]
-
\alpha(1-\alpha)
d_c^\alpha(F_X,G_Y)
\right\}.
\]

\begin{equation}
\begin{aligned}
d_{R-T}^{1/2}(F_X,G_Y)
&=
-\frac{2}{3}
\Bigg\{
\log
\iint
f_X^{1/2}(x_1,x_2)
g_Y^{1/2}(x_1,x_2)
\,dx_1\,dx_2
\\
&\qquad\qquad
+
\iint
f_X^{1/2}(x_1,x_2)
g_Y^{1/2}(x_1,x_2)
\,dx_1\,dx_2
-1
\Bigg\}
\\
&=
-\frac{2}{3}
\Bigg\{
\log
\iint
\sqrt{f_X(x_1,x_2)}
\sqrt{g_Y(x_1,x_2)}
\,dx_1\,dx_2
\\
&\qquad\qquad
+
\iint
\sqrt{f_X(x_1,x_2)}
\sqrt{g_Y(x_1,x_2)}
\,dx_1\,dx_2
-1
\Bigg\}.
\end{aligned}
\label{eq:half}
\end{equation}
The Bhattacharya divergence measure proposed by \cite{Bhattacharyya1943} is given by
\[
B(F_X,G_Y)
=
\iint
\sqrt{f_X(x_1,x_2)}
\sqrt{g_Y(x_1,x_2)}
\,dx_1\,dx_2.
\]
Therefore, equation \eqref{eq:half} can be written as
\begin{equation}
d_{R-T}^{1/2}(F_X,G_Y)
=
-\frac{2}{3}
\left[
\log B(F_X,G_Y)
+
B(F_X,G_Y)
-1
\right].
\label{eq:halfb}
\end{equation}
\cite{Hellinger1907} developed the divergence measure known as Hellinger divergence measure, given by
\begin{equation*}
\begin{aligned}
H^2(F_X,G_Y)
&=
\frac{1}{2}
\iint
\left(
\sqrt{f_X(x_1,x_2)}
-
\sqrt{g_Y(x_1,x_2)}
\right)^2
\,dx_1\,dx_2\\
&=
\frac{1}{2}
\iint
\left(
f_X(x_1,x_2)
-
2\sqrt{f_X(x_1,x_2)}
\sqrt{g_Y(x_1,x_2)}
+
g_Y(x_1,x_2)
\right)
\,dx_1\,dx_2\\
&=
\frac{1}{2}
\left[
1
-
2
\iint
\sqrt{f_X(x_1,x_2)}
\sqrt{g_Y(x_1,x_2)}
\,dx_1\,dx_2
+
1
\right]\\
&=
1
-
\iint
\sqrt{f_X(x_1,x_2)}
\sqrt{g_Y(x_1,x_2)}
\,dx_1\,dx_2.
\end{aligned}
\end{equation*}
Therefore,
$
B(F_X,G_Y)
=
1-H^2(F_X,G_Y).
$
Hence, equation \eqref{eq:halfb} becomes
\[
d_{R-T}^{1/2}(F_X,G_Y)
=
-\frac{2}{3}
\left[
\log
\left(
1-H^2(F_X,G_Y)
\right)
-
H^2(F_X,G_Y)
\right].
\]
Moreover,
$
d_{R-T}^{1/2}(F_X,G_Y)
=
-\frac{2}{3}
\left[
-d_B(F_X,G_Y)
+
\exp\left(
-d_B(F_X,G_Y)
\right)
-1
\right],
$
where
$
d_B(F_X,G_Y)
=
-\log B(F_X,G_Y)
$
is known as the Bhattacharya distance.

In the next example, the bivariate differential Rényi–Tsallis divergence measure of independent bivariate exponential distribution is compared with some
other divergence measures.

\begin{example}
Let \(X=(X_1,X_2)\) follow the independent bivariate exponential distribution with density
$
f(x_1,x_2)
=
\lambda_1\lambda_2
\exp(-\lambda_1x_1-\lambda_2x_2),
\ x_1,x_2\ge0,\ \lambda_1,\lambda_2>0.
$
Similarly, let
\(Y=(Y_1,Y_2)\) follow the independent bivariate exponential distribution with density
$
g(y_1,y_2)
=
\mu_1\mu_2
\exp(-\mu_1y_1-\mu_2y_2),
\ y_1,y_2\ge0,\ \mu_1,\mu_2>0.
$
We obtain the Rényi--Tsallis divergence measure as,

\begin{equation}
d_{R-T}^{\alpha}(F_X,G_Y)
=
\begin{cases}
\displaystyle
\frac{1}{\alpha-\alpha^{-1}}
\Bigg[
\log
\Bigg\{
\prod_{i=1}^{2}
\left(
\frac{
\lambda_i^{\alpha}\mu_i^{1-\alpha}
}{
\alpha\lambda_i+(1-\alpha)\mu_i
}
\right)
\Bigg\}
\\[0.4cm]
\displaystyle\qquad
+
\prod_{i=1}^{2}
\left(
\frac{
\lambda_i^{\alpha}\mu_i^{1-\alpha}
}{
\alpha\lambda_i+(1-\alpha)\mu_i
}
\right)
-1
\Bigg],
\qquad \alpha > 0,\alpha\neq1
\\[0.8cm]

\displaystyle
\sum_{i=1}^{2}
\left[
\log\left(\frac{\lambda_i}{\mu_i}\right)
+
\frac{\mu_i-\lambda_i}{\lambda_i}
\right],
\qquad \alpha=1.
\end{cases}
\end{equation}
\end{example}

\begin{figure}[htbp]
\centering

\subfloat[]{%
\resizebox*{0.32\textwidth}{!}{\includegraphics{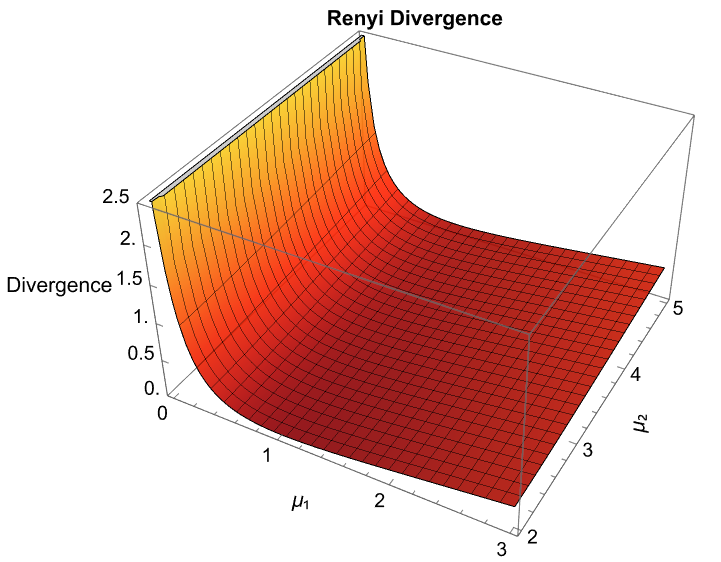}}}
\hspace{2pt}
\subfloat[]{%
\resizebox*{0.32\textwidth}{!}{\includegraphics{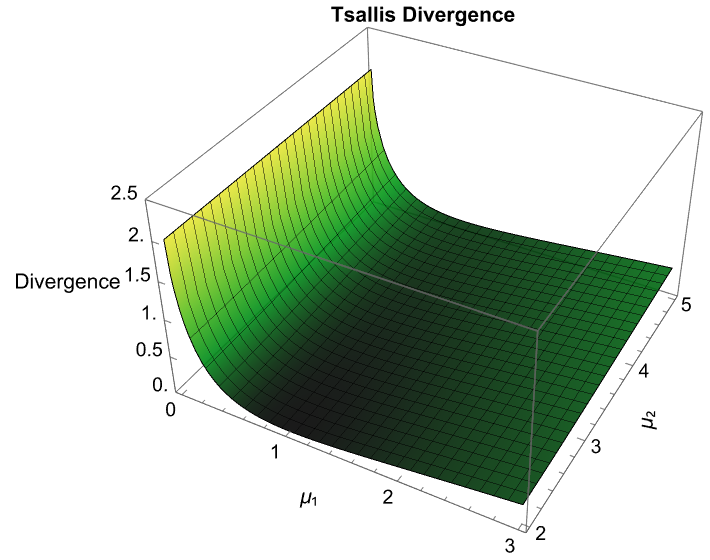}}}
\hspace{2pt}
\subfloat[]{%
\resizebox*{0.32\textwidth}{!}{\includegraphics{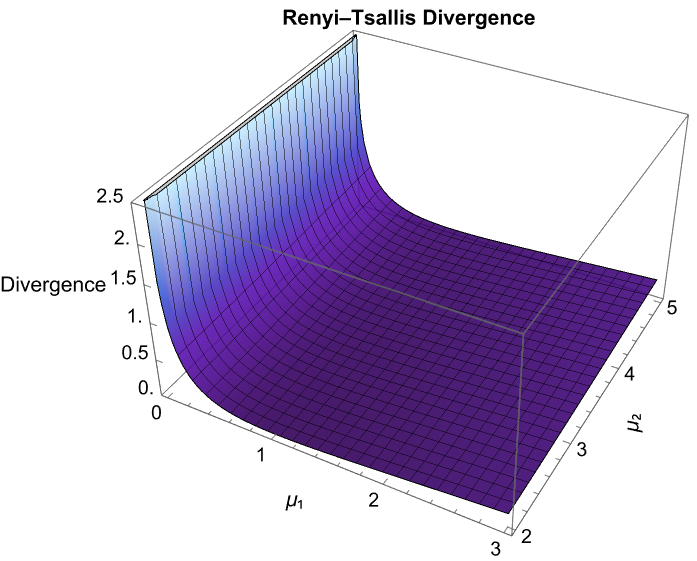}}}

\caption{Plots of Rényi, Tsallis, and Rényi--Tsallis divergence measures for independent bivariate exponential distributions with $\alpha = 0.5$, $\lambda_1 = 1$, $\lambda_2 = 2$, $\mu_1 \in (0,3]$, and $\mu_2 \in [2,5]$.}
\label{fig:divergence-plots}
\end{figure}

Figure~\ref{fig:divergence-plots} shows the Rényi, Tsallis and Rényi--Tsallis divergence measures for independent bivariate exponential distribution with $\alpha = 0.5$, $\lambda_1 = 1$, $\lambda_2 = 2$,
$\mu_1 \in [0,3]$, and $\mu_2 \in [2,5]$. All three measures increase with increasing $\mu_2$ and decrease with 
increasing $\mu_1$. They take their largest values for small $\mu_1$ and large 
$\mu_2$. They 
all exhibit the same overall monotonic trend; however, the Rényi divergence is 
relatively more sensitive to changes in the parameters, the Tsallis divergence 
changes more smoothly, and the proposed Rényi--Tsallis divergence maintains 
the common sensitivity pattern while providing a balanced measure of the 
distributional difference.
And Table~\ref{tab:divergence_measures} gives the different divergence measures of independent bivariate exponential distribution. Following proposition gives some useful bound for the proposed divergence measure.

\begin{table}
\centering
\tbl{Some divergence measures for independent bivariate exponential distributions}
{\begin{tabular}{p{5cm} p{5cm}}
\toprule
\textbf{Divergence Measure} & \textbf{Expression} \\
\midrule

Kullback--Leibler Divergence
&
$\displaystyle
\sum_{i=1}^{2}
\left[
\log\left(\frac{\lambda_i}{\mu_i}\right)
+
\frac{\mu_i-\lambda_i}{\lambda_i}
\right]
$
\\[0.4cm]

R\'enyi Divergence
&
$\displaystyle
\frac{1}{\alpha-1}
\log
\left[
\prod_{i=1}^{2}
\left(
\frac{\lambda_i^{\alpha}\mu_i^{1-\alpha}}
{\alpha\lambda_i+(1-\alpha)\mu_i}
\right)
\right]
$
\\[0.4cm]

Tsallis Divergence
&
$\displaystyle
\frac{1}{\alpha-1}
\left[
\prod_{i=1}^{2}
\left(
\frac{\lambda_i^{\alpha}\mu_i^{1-\alpha}}
{\alpha\lambda_i+(1-\alpha)\mu_i}
\right)
-1
\right]
$
\\[0.4cm]

$\alpha$-Divergence
&
$\displaystyle
\frac{1}{\alpha(\alpha-1)}
\left[
1-
\prod_{i=1}^{2}
\left(
\frac{\lambda_i^{\alpha}\mu_i^{1-\alpha}}
{\alpha\lambda_i+(1-\alpha)\mu_i}
\right)
\right]
$
\\[0.4cm]

Hellinger Distance
&
$\displaystyle
1-
\prod_{i=1}^{2}
\left(
\frac{2\sqrt{\lambda_i\mu_i}}
{\lambda_i+\mu_i}
\right)
$
\\[0.4cm]

Bhattacharyya Distance
&
$\displaystyle
-\log
\left[
\prod_{i=1}^{2}
\left(
\frac{2\sqrt{\lambda_i\mu_i}}
{\lambda_i+\mu_i}
\right)
\right]
$
\\
\bottomrule
\end{tabular}}
\label{tab:divergence_measures}
\end{table}

\begin{example}
Sensitivity of divergence measures under dependent bivariate models.\\
The above example shows the behaviour of the Rényi, Tsallis, and Rényi–Tsallis divergence measures when the parameters of independent bivariate exponential distributions vary. In this example, we study their sensitivity to changes in dependence under two different bivariate models: the Marshall–Olkin model and a Gumbel-type bivariate exponential model. 

As a common reference distribution, consider the independent bivariate exponential distribution with unit marginal rates, whose density is
$
q(x_1,x_2)=e^{-x_1-x_2},\qquad x_1,x_2\geq 0.$
For a dependent model with density \(p\), we evaluate the Rényi, Tsallis, and Rényi–Tsallis divergence measures using equations \eqref{eq:Rényi}, \eqref{eq:Tsallis} and \eqref{eq:RT}.\\
\textbf{Marshall–Olkin bivariate exponential model}\\
Let \((X_1,X_2)\) follow a Marshall–Olkin bivariate exponential distribution with survival function
\[
\overline F_c(x_1,x_2)
=\exp\left[-(1-c)x_1-(1-c)x_2
-c\max(x_1,x_2)\right],
\]
where \(x_1,x_2\geq0\) and \(0\leq c<1\). The marginal distributions are standard exponential for every value of \(c\). The parameter \(c\) controls the common-shock component, and \(c=0\) corresponds to independence.
The absolutely continuous part of the joint distribution has density
\[
f_c(x_1,x_2)=
\begin{cases}
(1-c)e^{-(1-c)x_1-x_2},&x_1<x_2,\\[2mm]
(1-c)e^{-x_1-(1-c)x_2},&x_2<x_1.
\end{cases}
\]
For \(c>0\), the distribution also has a singular component on the diagonal \(x_1=x_2\), with density \(ce^{-(2-c)t}\) with respect to the diagonal parameter \(t\geq0\). Its total probability mass is \(c/(2-c)\).

Using a common dominating measure that accounts for both the continuous and singular components, the power integral relative to the independent reference distribution is
\[
I_\alpha^{MO}(c)
=\frac{2(1-c)^\alpha}{2-\alpha c},
\qquad 0<\alpha<1.
\]
Consequently, the three divergence measures can be obtained by substituting \(I_\alpha^{MO}(c)\) into their respective definitions. The diagonal component contributes no term to this power integral because the reference distribution assigns zero probability to the diagonal.\\
The Kendall rank correlation coefficient for this model is
\[
\tau_{MO}(c)=\frac{c}{2-c}.
\]
Writing \(\rho=|\tau|\), the common-shock parameter can therefore be expressed as
\[
c=\frac{2\rho}{1+\rho}.
\]
This parameterization allows the sensitivity of the divergence measures to be evaluated with respect to the magnitude of Kendall's tau.\\
\textbf{Gumbel-type bivariate exponential model}\\
For comparison, consider a Gumbel-type bivariate exponential model with survival function
\[
\overline G_\theta(x_1,x_2)
=\exp(-x_1-x_2-\theta x_1x_2),
\qquad 0\leq\theta\leq1.
\]
Both marginal distributions are standard exponential, and \(\theta=0\) corresponds to independence. Its joint density is
\[
g_\theta(x_1,x_2)
=e^{-x_1-x_2-\theta x_1x_2}
\left[(1+\theta x_1)(1+\theta x_2)-\theta\right].
\]
The power integral relative to the same independent reference distribution is
\[
I_\alpha^G(\theta)
=\int_0^\infty\int_0^\infty
e^{-(x_1+x_2)-\alpha\theta x_1x_2}
\times
\left[(1+\theta x_1)(1+\theta x_2)-\theta\right]^\alpha
\,dx_2\,dx_1.
\]
Unlike the corresponding integral for the Marshall–Olkin model, this expression does not simplify to the same closed form and is evaluated numerically using Gauss–Legendre quadrature.
The Kendall rank correlation coefficient is evaluated numerically from
\[
\tau_G(\theta)
=4\int_0^\infty\int_0^\infty
G_\theta(x_1,x_2)g_\theta(x_1,x_2)
\,dx_2\,dx_1-1,
\]
where
\[
G_\theta(x_1,x_2)
=1-e^{-x_1}-e^{-x_2}
+e^{-x_1-x_2-\theta x_1x_2}
\]
is the joint distribution function. In the parameter range considered, this coefficient is non-positive. We therefore compare the two models using the magnitude \(\rho=|\tau|\), matching the value of \(\rho\) in the Marshall–Olkin model to the corresponding value in the Gumbel-type model.\\
\textbf{Sensitivity with respect to dependence}\\
Let \(D^\alpha\) denote any one of the three divergence measures. For each model, define its absolute sensitivity with respect to the magnitude of Kendall's tau by
\[
S_D^\alpha(\rho)
=\left|\frac{dD^\alpha}{d\rho}\right|.
\]
For the Marshall–Olkin model, the chain rule gives
\[
S_{D,MO}^\alpha(\rho)
=\left|\frac{\partial D_{MO}^\alpha}{\partial c}\right|
\frac{2}{(1+\rho)^2}.
\]
For the Gumbel-type model, the corresponding sensitivity is
\[
S_{D,G}^\alpha(\rho)
=\left|
\frac{dD_G^\alpha/d\theta}
{d|\tau_G(\theta)|/d\theta}
\right|_{\theta=\theta(\rho)},
\]
where \(\theta(\rho)\) is obtained by numerically solving
\[
|\tau_G(\theta)|=\rho.
\]

We evaluate these sensitivities for \(\alpha=0.25,0.50,0.75,\) and \(0.90\) over the common range of dependence magnitudes attainable by both models. This comparison places the results on a common dependence scale instead of comparing derivatives with respect to the model's differently scaled parameters.

\begin{figure}[htbp]
\centering
\includegraphics[width=0.5\textwidth]
{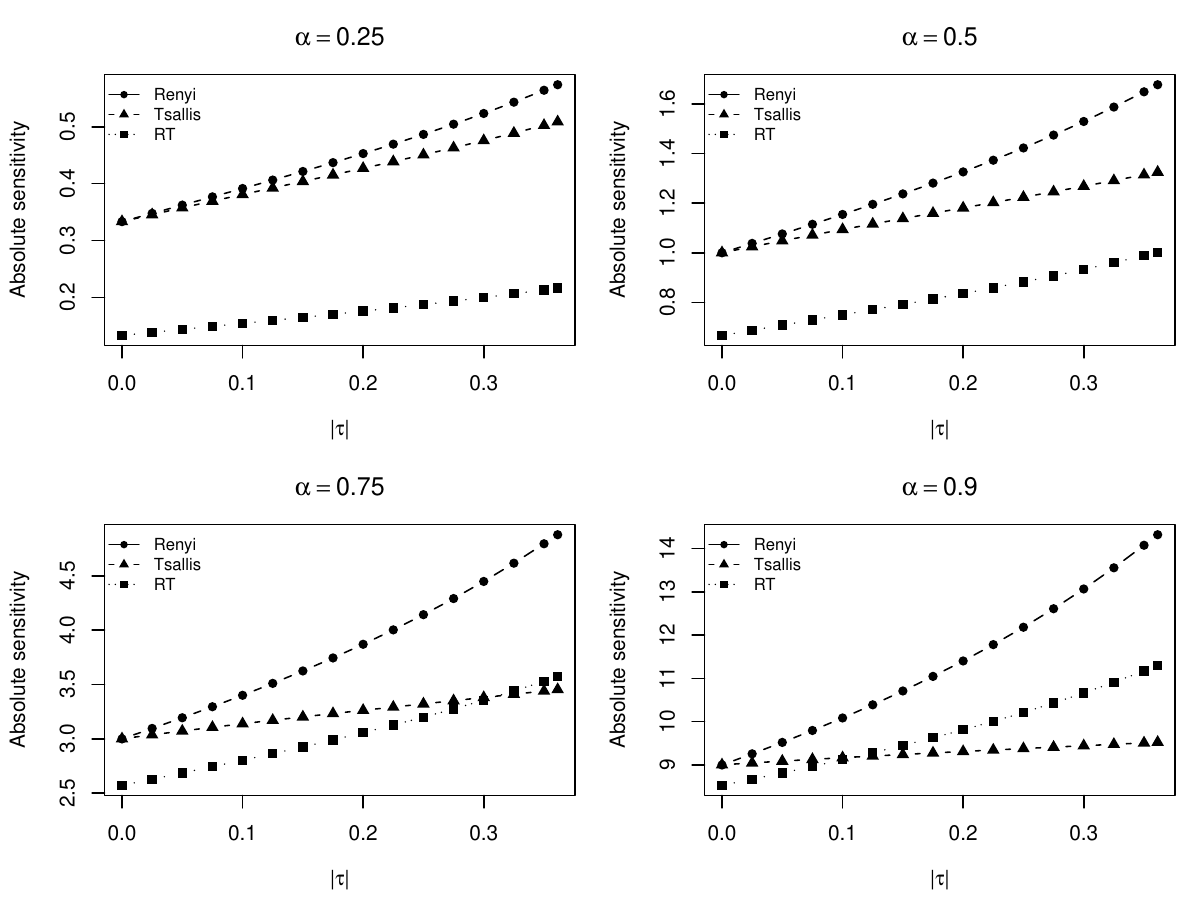}
\caption{Absolute sensitivity of the Rényi, Tsallis, and
Rényi--Tsallis divergence measures with respect to the magnitude
of Kendall's tau under the Marshall--Olkin bivariate exponential
model, for $\alpha=0.25,0.50,0.75,$ and $0.90$.}
\label{fig:MO-sensitivity}
\end{figure}

\begin{figure}[htbp]
\centering
\includegraphics[width=0.5\textwidth]
{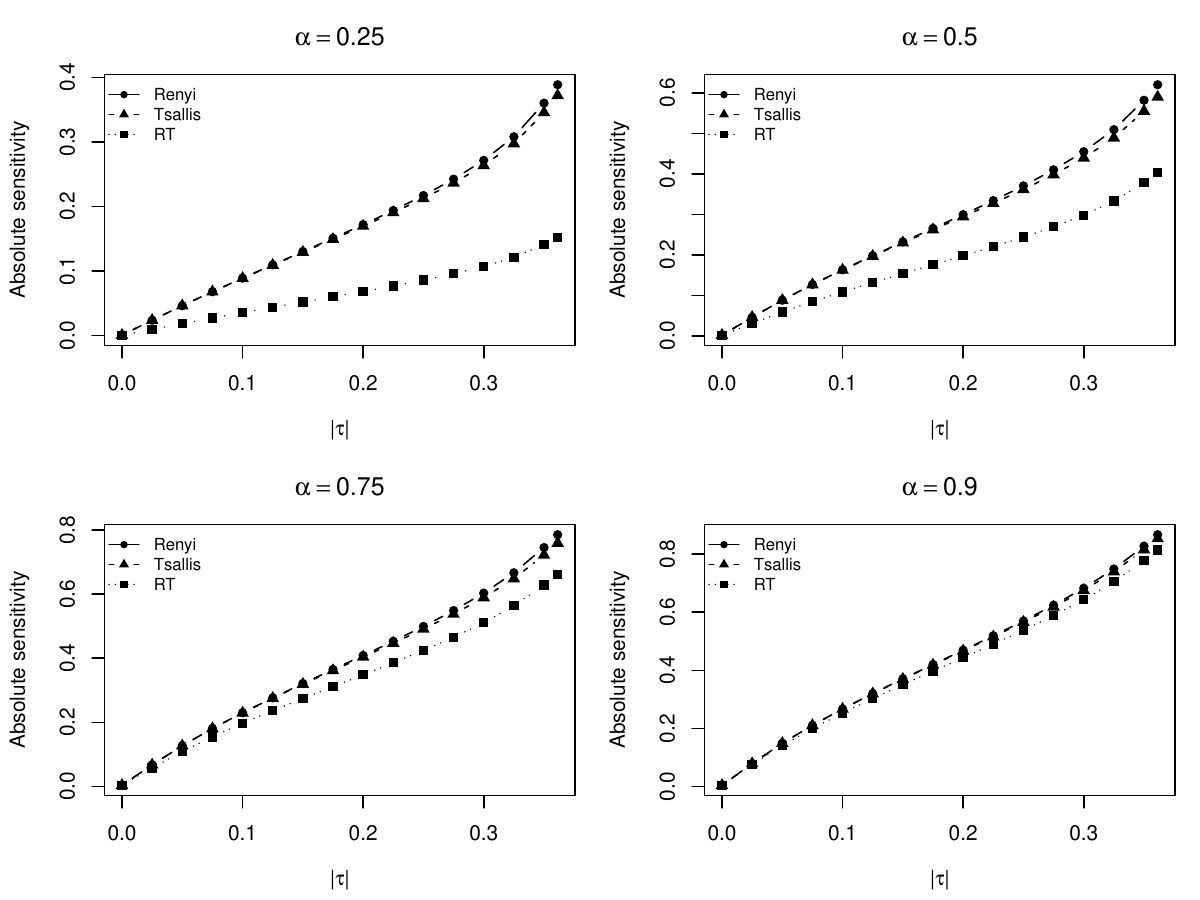}
\caption{Absolute sensitivity of the Rényi, Tsallis, and
Rényi--Tsallis divergence measures with respect to the magnitude
of Kendall's tau under the Gumbel-type bivariate exponential
model, for $\alpha=0.25,0.50,0.75,$ and $0.90$.}
\label{fig:Gumbel-sensitivity}
\end{figure}

\begin{figure}[htbp]
\centering
\includegraphics[width=0.5\textwidth]
{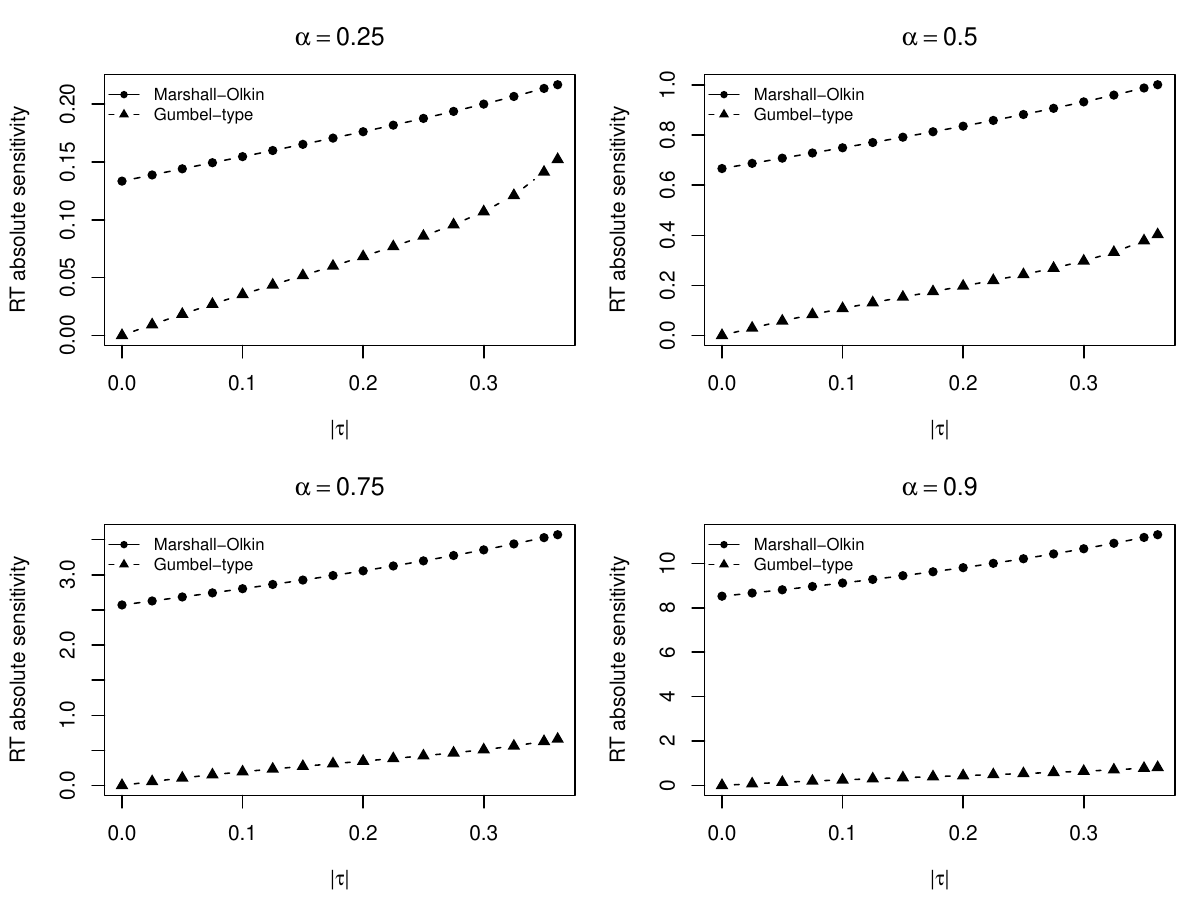}
\caption{Comparison of the absolute sensitivity of the
Rényi--Tsallis divergence under the Marshall--Olkin and
Gumbel-type bivariate exponential models at matched values
of the magnitude of Kendall's tau.}
\label{fig:RT-model-comparison}
\end{figure}

\begin{table}[htbp]
\centering
\caption{Mean absolute sensitivity of the divergence measures
over the matched dependence grid.}
\label{tab:mean-sensitivity}
{\scriptsize
\setlength{\tabcolsep}{4pt}
\renewcommand{\arraystretch}{0.9}
\begin{tabular}{ccccc}
\toprule
$\alpha$ & Divergence & MO & Gumbel-type & Ratio (MO/Gumbel)\\
\midrule
0.25 & R\'enyi   & 0.4497 & 0.1732 & 2.60\\
0.25 & Tsallis   & 0.4222 & 0.1689 & 2.50\\
0.25 & R\'enyi--Tsallis & 0.1744 & 0.0684 & 2.55\\
\midrule
0.50 & R\'enyi   & 1.3211 & 0.2938 & 4.50\\
0.50 & Tsallis   & 1.1690 & 0.2855 & 4.10\\
0.50 & R\'enyi--Tsallis & 0.8300 & 0.1931 & 4.30\\
\midrule
0.75 & R\'enyi   & 3.8693 & 0.3926 & 9.86\\
0.75 & Tsallis   & 3.2417 & 0.3850 & 8.42\\
0.75 & R\'enyi--Tsallis & 3.0476 & 0.3333 & 9.14\\
\midrule
0.90 & R\'enyi   & 11.4226 & 0.4457 & 25.63\\
0.90 & Tsallis   & 9.2824 & 0.4418 & 21.01\\
0.90 & R\'enyi--Tsallis & 9.8076 & 0.4204 & 23.33\\
\bottomrule
\end{tabular}}
\end{table}
Table~\ref{tab:mean-sensitivity} shows that the mean absolute sensitivities of all three divergence measures are higher under the Marshall--Olkin model than under the Gumbel-type model for the different values of $\alpha$. The difference increases with $\alpha$; at $\alpha=0.90$, the ratios are 25.63, 21.01, and 23.33 for the Rényi, Tsallis, and Rényi--Tsallis divergences, respectively. However, the ordering of relative sensitivity is different for different models and orders, showing that sensitivity depends on both the divergence measure and the dependence structure. Marshall--Olkin model has a singular diagonal component whereas the Gumbel-type model is absolutely continuous. These differences should not be attributed only to the divergence measures. Numerical accuracy was tested through quadrature convergence, verification of the Rényi--Tsallis identity, and matching of the dependence magnitudes using $|\tau|$.
\end{example}

\begin{proposition}
Let $X=(x_1,x_2)$ and $Y=(y_1,y_2)$ be two non-negative absolutely continuous random vectors with density functions $f_X(x_1,x_2)$ and $g_Y(x_1,x_2)$ respectively. Since for all
$
f_X(x_1,x_2),g_Y(x_1,x_2)\ge0,
$
we have
\[
\log
\left(
\iint
f_X^\alpha(x_1,x_2)
g_Y^{\,1-\alpha}(x_1,x_2)
\,dx_1\,dx_2
\right)
\le
\left(
\iint
f_X^\alpha(x_1,x_2)
g_Y^{\,1-\alpha}(x_1,x_2)
\,dx_1\,dx_2
\right)
-1,
\ \alpha>0.
\]
Also,
$
d_{R-T}^{\alpha}(F_X,G_Y)
=
\frac{\alpha}{\alpha+1}
\left(
d_R^\alpha+d_T^\alpha
\right).
$
If $0<\alpha<1$, then
$
d_R^\alpha\ge d_T^\alpha.
$
Therefore,
$
2d_T^\alpha
\le
d_R^\alpha+d_T^\alpha
\le
2d_R^\alpha.
$
Hence,
\[
\frac{2\alpha}{\alpha+1}
d_T^\alpha
\le
d_{R-T}^\alpha
\le
\frac{2\alpha}{\alpha+1}
d_R^\alpha.
\]
If $\alpha>1$, then
$
d_R^\alpha\le d_T^\alpha.
$
Therefore,
\[
\frac{2\alpha}{\alpha+1}
d_R^\alpha
\le
d_{R-T}^\alpha
\le
\frac{2\alpha}{\alpha+1}
d_T^\alpha.
\]
In addition, \cite{Dragomir2000} obtained the bound
\[d_{KL}(F_X,G_Y)
\le
\log\left(
d_{\chi^2}(F_X,G_Y)+1
\right)\]
where $d_{\chi^2}(F_X,G_Y)$ is known as the $\chi^2$-distance (\cite{BenTal1991}).
Hence,
\begin{equation}
\left\{
\begin{alignedat}{2}
&\frac{2\alpha}{\alpha+1}
\left[d_T^\alpha(F_X,G_Y)\right]
\le
d_{R-T}^\alpha(F_X,G_Y)
\le
\frac{2\alpha}{\alpha+1}
\left[d_R^\alpha(F_X,G_Y)\right],
&\qquad& 0<\alpha<1,
\\[0.4cm]
&
d_{R-T}^\alpha(F_X,G_Y)
\le
\log\left[
d_{\chi^2}(F_X,G_Y)+1
\right],
&\qquad& \alpha=1,
\\[0.4cm]
&\frac{2\alpha}{\alpha+1}
\left[d_R^\alpha(F_X,G_Y)\right]
\le
d_{R-T}^\alpha(F_X,G_Y)
\le
\frac{2\alpha}{\alpha+1}
\left[d_T^\alpha(F_X,G_Y)\right],
&\qquad& \alpha>1.
\end{alignedat}
\right.
\label{eq:bound}
\end{equation}
\end{proposition}
Equation (\ref{eq:bound}) shows that the bivariate Rényi--Tsallis divergence lies between scaled versions of bivariate Tsallis divergence and bivariate Rényi divergence. The scale factor is $\frac{2\alpha}{\alpha+1}$ which lies between 0 and 2 for all $\alpha>0$. This bound shows how the interval containing the proposed divergence measure varies with the order parameter $\alpha$.  Also for $\alpha=1$, $\frac{2\alpha}{\alpha+1}$ becomes 1 and all these three divergences will be same and equal to KL divergence.
\section{Rényi-Tsallis Divergence for Bivariate Order Statistics}
\label{sec:order-statistics}

Let $X_1,X_2,\ldots,X_n$ and $Y_1,Y_2,\ldots,Y_n$ be two independent random samples of size $n$ from absolutely continuous distributions with cdfs $F_X$ and $G_Y$, and corresponding pdfs $f_X$ and $g_Y$, respectively. Let $X_{(i)},X_{(j)}$ and $Y_{(i)},Y_{(j)}$ denote the $i^{th}$ and $j^{th}$ order statistics $(i<j)$ from the samples associated with $F_X$ and $G_Y$, respectively.

For $0<\alpha<\infty$, $\alpha\neq 1$, the Rényi--Tsallis divergence between the bivariate order statistic distributions $F_{ij}$ and $G_{ij}$ is defined as
\begin{equation}
\begin{aligned}
d_{R\text{-}T}^{\alpha}(F_X,G_Y)
&=
\frac{1}{\alpha-\alpha^{-1}}
\Bigg[
\log
\left(
\int_{0}^{\infty}\int_{x_i}^{\infty}
f_{ij}^{\alpha}(x_i,x_j)\,
g_{ij}^{\,1-\alpha}(x_i,x_j)
\,dx_j\,dx_i
\right)
\\
&\qquad\qquad
+
\int_{0}^{\infty}\int_{x_i}^{\infty}
f_{ij}^{\alpha}(x_i,x_j)\,
g_{ij}^{\,1-\alpha}(x_i,x_j)
\,dx_j\,dx_i
-1
\Bigg].
\end{aligned}
\label{eq:RT_ordered}
\end{equation}
The joint pdf of the pair $(X_{(i)},X_{(j)})$ is given by
\[
\begin{aligned}
f_{ij}(x_i,x_j)
&=
\frac{n!}{(i-1)!(j-i-1)!(n-j)!}
[F_X(x_i)]^{i-1}
[F_X(x_j)-F_X(x_i)]^{j-i-1}
\\
&\qquad\times
[1-F_X(x_j)]^{n-j}
\,f_X(x_i)f_X(x_j),
\qquad x_i<x_j .
\end{aligned}
\]
Similarly, the joint pdf of the pair $(Y_{(i)},Y_{(j)})$ is
\[
\begin{aligned}
g_{ij}(y_i,y_j)
&=
\frac{n!}{(i-1)!(j-i-1)!(n-j)!}
[G_Y(y_i)]^{i-1}
[G_Y(y_j)-G_Y(y_i)]^{j-i-1}
\\
&\qquad\times
[1-G_Y(y_j)]^{n-j}
\,g_Y(y_i)g_Y(y_j),
\qquad y_i<y_j .
\end{aligned}
\]
The differential Rényi--Tsallis divergence between the bivariate order statistics of two samples admits a joint representation in terms of the Rényi divergence and the Tsallis divergence based on their joint probability density functions. That is,
\begin{equation}
d_{R\text{-}T}^{\alpha}(F_X,G_Y)
=
\frac{\alpha}{\alpha+1}
\left[
d_{R}^{\alpha}(F_X,G_Y)
+
d_{T}^{\alpha}(F_X,G_Y)
\right],
\end{equation}
where $d_{R}^{\alpha}(F_X,G_Y)$ and $d_{T}^{\alpha}(F_X,G_Y)$ denote the Rényi divergence and Tsallis divergence, respectively, between the joint densities of the bivariate order statistics given by,
\[
d_{R}^{\alpha}(F_X,G_Y)
=
\frac{1}{\alpha-1}
\log
\left[
\int_{0}^{\infty}\int_{x_i}^{\infty}
f_{ij}^{\alpha}(x_i,x_j)
\,g_{ij}^{\,1-\alpha}(x_i,x_j)
\,dx_j\,dx_i
\right],
\]

\[
d_{T}^{\alpha}(F_X,G_Y)
=
\frac{1}{\alpha-1}
\left[
\int_{0}^{\infty}\int_{x_i}^{\infty}
f_{ij}^{\alpha}(x_i,x_j)
\,g_{ij}^{\,1-\alpha}(x_i,x_j)
\,dx_j\,dx_i
-1
\right].
\]

\subsection*{Rényi--Tsallis Divergence for Bivariate Extreme Order Statistics}

Let $X_1,X_2,\ldots,X_n$ be a random sample from a continuous distribution with cumulative distribution function $F_X(x)$ and probability density function $f_X(x)$. Similarly, let $Y_1,Y_2,\ldots,Y_n$ be a random sample from another continuous distribution with cumulative distribution function $G_Y(x)$ and probability density function $g_Y(x)$.

Let $X_{1:n}$ and $X_{n:n}$ denote the minimum and maximum order statistics, respectively. Then, for $x_1<x_n$, the joint probability density function of the bivariate extreme order statistic vector $(X_{1:n},X_{n:n})$ is given by
\[
f_{1,n}(x_1,x_n)
=
n(n-1)
\left[F_X(x_n)-F_X(x_1)\right]^{n-2}
f_X(x_1)f_X(x_n).
\]
Similarly,
\[
g_{1,n}(x_1,x_n)
=
n(n-1)
\left[G_Y(x_n)-G_Y(x_1)\right]^{n-2}
g_Y(x_1)g_Y(x_n).
\]
Then, the Rényi--Tsallis divergence of order $\alpha>0$, $\alpha\neq 1$, between the bivariate extreme order statistic distributions corresponding to $F_X$ and $G_Y$ is defined as
\begin{equation}
d_{R\text{-}T}^{\alpha}(F_X,G_Y)
=
\frac{1}{\alpha-\alpha^{-1}}
\left[
\log I_{\alpha}
+
I_{\alpha}
-1
\right],
\end{equation}
where
\begin{equation}
\begin{aligned}
I_{\alpha}
&=
\int_{-\infty}^{\infty}
\int_{x_1}^{\infty}
f_{1,n}^{\alpha}(x_1,x_n)\,
g_{1,n}^{\,1-\alpha}(x_1,x_n)
\,dx_n\,dx_1 \\
&=
n(n-1)
\int_{-\infty}^{\infty}
\int_{x_1}^{\infty}
[F_X(x_n)-F_X(x_1)]^{\alpha(n-2)}
[G_Y(x_n)-G_Y(x_1)]^{(1-\alpha)(n-2)}
\\
&\qquad\qquad\times
f_X^{\alpha}(x_1)\,
g_Y^{\,1-\alpha}(x_1)\,
f_X^{\alpha}(x_n)\,
g_Y^{\,1-\alpha}(x_n)
\,dx_n\,dx_1 .
\end{aligned}
\end{equation}

\begin{example}: \textbf{Exponential Distribution}\\
Let
$
X_1,\ldots,X_n \stackrel{\text{i.i.d.}}{\sim} \mathrm{Exp}(\lambda_1), \
Y_1,\ldots,Y_n \stackrel{\text{i.i.d.}}{\sim} \mathrm{Exp}(\lambda_2),
$
where $\lambda_1>0$ and $\lambda_2>0$. The corresponding density and distribution functions are
\[
f(x)=\lambda_1 e^{-\lambda_1 x},
\qquad
F_X(x)=1-e^{-\lambda_1 x},
\qquad x>0,
\]
and
\[
g(x)=\lambda_2 e^{-\lambda_2 x},
\qquad
G_Y(x)=1-e^{-\lambda_2 x},
\qquad x>0.
\]

For the bivariate extreme order statistics $(X_{1:n},X_{n:n})$ and $(Y_{1:n},Y_{n:n})$, the joint densities are
\[
\begin{aligned}
f_{1,n}(x_1,x_n)
&=
n(n-1)
\left(
e^{-\lambda_1 x_1}
-
e^{-\lambda_1 x_n}
\right)^{n-2}
\lambda_1^2
e^{-\lambda_1(x_1+x_n)},
\end{aligned}
\]
and
\[
\begin{aligned}
g_{1,n}(x_1,x_n)
&=
n(n-1)
\left(
e^{-\lambda_2 x_1}
-
e^{-\lambda_2 x_n}
\right)^{n-2}
\lambda_2^2
e^{-\lambda_2(x_1+x_n)}.
\end{aligned}
\]

Then the Rényi--Tsallis divergence of order $\alpha>0$, $\alpha\neq 1$ between these two bivariate extreme order statistics is given by,
\begin{equation}
d_{R\text{-}T}^{\alpha}(F_X,G_Y)
=
\frac{1}{\alpha-\alpha^{-1}}
\left[
\log(I_\alpha)
+
I_\alpha
-1
\right],
\end{equation}
where
\begin{equation}
I_\alpha
=
\int_0^\infty
\int_{x_1}^{\infty}
f_{1,n}^{\alpha}(x_1,x_n)\,
g_{1,n}^{\,1-\alpha}(x_1,x_n)
\,dx_n\,dx_1 .
\end{equation}
Substituting the corresponding joint densities, we obtain
\begin{equation}
\begin{aligned}
I_\alpha
&=
n(n-1)
\lambda_1^{2\alpha}
\lambda_2^{2(1-\alpha)}
\int_0^\infty
\int_{x_1}^{\infty}
\left(
e^{-\lambda_1 x_1}
-
e^{-\lambda_1 x_n}
\right)^{\alpha(n-2)}
\\
&\quad\times
\left(
e^{-\lambda_2 x_1}
-
e^{-\lambda_2 x_n}
\right)^{(1-\alpha)(n-2)}
e^{-[\alpha\lambda_1+(1-\alpha)\lambda_2](x_1+x_n)}
\,dx_n\,dx_1 .
\end{aligned}
\end{equation}
\end{example}

\begin{example}
\textbf{Weibull Distribution}\\

Let
$
X_1,\ldots,X_n \stackrel{\text{i.i.d.}}{\sim} \mathrm{Weibull}(\lambda_1,\beta_1),
\
Y_1,\ldots,Y_n \stackrel{\text{i.i.d.}}{\sim} \mathrm{Weibull}(\lambda_2,\beta_2),
$
where $\lambda_1,\lambda_2,\beta_1,\beta_2>0$. The corresponding density and distribution functions are
\[
f(x)
=
\beta_1\lambda_1^{\beta_1}
x^{\beta_1-1}
e^{-(\lambda_1 x)^{\beta_1}},
\qquad
F_X(x)
=
1-e^{-(\lambda_1 x)^{\beta_1}},
\qquad x>0,
\]
and
\[
g(x)
=
\beta_2\lambda_2^{\beta_2}
x^{\beta_2-1}
e^{-(\lambda_2 x)^{\beta_2}},
\qquad
G_Y(x)
=
1-e^{-(\lambda_2 x)^{\beta_2}},
\qquad x>0.
\]
For the bivariate extreme order statistics $(X_{1:n},X_{n:n})$ and $(Y_{1:n},Y_{n:n})$, the joint densities are
\[
\begin{aligned}
f_{1,n}(x_1,x_n)
&=
n(n-1)
\left(
e^{-(\lambda_1 x_1)^{\beta_1}}
-
e^{-(\lambda_1 x_n)^{\beta_1}}
\right)^{n-2}
\\
&\quad\times
\beta_1^2\lambda_1^{2\beta_1}
x_1^{\beta_1-1}
x_n^{\beta_1-1}
e^{-(\lambda_1 x_1)^{\beta_1}-(\lambda_1 x_n)^{\beta_1}},
\end{aligned}
\]
and
\[
\begin{aligned}
g_{1,n}(x_1,x_n)
&=
n(n-1)
\left(
e^{-(\lambda_2 x_1)^{\beta_2}}
-
e^{-(\lambda_2 x_n)^{\beta_2}}
\right)^{n-2}
\\
&\quad\times
\beta_2^2\lambda_2^{2\beta_2}
x_1^{\beta_2-1}
x_n^{\beta_2-1}
e^{-(\lambda_2 x_1)^{\beta_2}-(\lambda_2 x_n)^{\beta_2}}.
\end{aligned}
\]
Then the Rényi--Tsallis divergence of order $\alpha>0$, $\alpha\neq 1$ between these two bivariate extreme order statistics is given by
\begin{equation}
d_{R\text{-}T}^{\alpha}(F_X,G_Y)
=
\frac{1}{\alpha-\alpha^{-1}}
\left[
\log I_\alpha
+
I_\alpha
-1
\right],
\end{equation}
where
\begin{equation}
I_\alpha
=
\int_0^\infty
\int_{x_1}^{\infty}
f_{1,n}^{\alpha}(x_1,x_n)\,
g_{1,n}^{\,1-\alpha}(x_1,x_n)
\,dx_n\,dx_1 .
\end{equation}
Substituting the corresponding joint densities, we obtain
\begin{equation}
\begin{aligned}
I_\alpha
&=
n(n-1)
\beta_1^{2\alpha}
\beta_2^{2(1-\alpha)}
\lambda_1^{2\alpha\beta_1}
\lambda_2^{2(1-\alpha)\beta_2}
\int_0^\infty
\int_{x_1}^{\infty}
\left(
e^{-(\lambda_1 x_1)^{\beta_1}}
-
e^{-(\lambda_1 x_n)^{\beta_1}}
\right)^{\alpha(n-2)}
\\
&\quad\times
\left(
e^{-(\lambda_2 x_1)^{\beta_2}}
-
e^{-(\lambda_2 x_n)^{\beta_2}}
\right)^{(1-\alpha)(n-2)}
(x_1 x_n)^{\alpha(\beta_1-1)+(1-\alpha)(\beta_2-1)}
\\
&\quad\times
e^{-\alpha[(\lambda_1 x_1)^{\beta_1}+(\lambda_1 x_n)^{\beta_1}]}
e^{-(1-\alpha)[(\lambda_2 x_1)^{\beta_2}+(\lambda_2 x_n)^{\beta_2}]}
\,dx_n\,dx_1 .
\end{aligned}
\end{equation}
\end{example}

\begin{example}
\textbf{Pareto Type-I Distribution}\\

Let
$
X_1,\ldots,X_n \stackrel{\text{i.i.d.}}{\sim} \mathrm{Pareto}(\lambda_1,\beta_1),
\
Y_1,\ldots,Y_n \stackrel{\text{i.i.d.}}{\sim} \mathrm{Pareto}(\lambda_2,\beta_2),
$
where $\lambda_1,\lambda_2,\beta_1,\beta_2>0$. The corresponding density and distribution functions are
\[
f(x)
=
\beta_1\lambda_1^{\beta_1}x^{-(\beta_1+1)},
\qquad
F_X(x)
=
1-\left(\frac{\lambda_1}{x}\right)^{\beta_1},
\qquad x\ge \lambda_1,
\]
and
\[
g(x)
=
\beta_2\lambda_2^{\beta_2}x^{-(\beta_2+1)},
\qquad
G_Y(x)
=
1-\left(\frac{\lambda_2}{x}\right)^{\beta_2},
\qquad x\ge \lambda_2.
\]

For the bivariate extreme order statistics $(X_{1:n},X_{n:n})$ and $(Y_{1:n},Y_{n:n})$, the joint densities are
\[
\begin{aligned}
f_{1,n}(x_1,x_n)
=
n(n-1)
\left[
\left(\frac{\lambda_1}{x_1}\right)^{\beta_1}
-
\left(\frac{\lambda_1}{x_n}\right)^{\beta_1}
\right]^{n-2}
\beta_1^2\lambda_1^{2\beta_1}
(x_1x_n)^{-(\beta_1+1)},
\end{aligned}
\]
and
\[
\begin{aligned}
g_{1,n}(x_1,x_n)
=
n(n-1)
\left[
\left(\frac{\lambda_2}{x_1}\right)^{\beta_2}
-
\left(\frac{\lambda_2}{x_n}\right)^{\beta_2}
\right]^{n-2}
\beta_2^2\lambda_2^{2\beta_2}
(x_1x_n)^{-(\beta_2+1)}.
\end{aligned}
\]
Then the Rényi--Tsallis divergence of order $\alpha>0$, $\alpha\neq 1$ between these two bivariate extreme order statistics is given by
\begin{equation}
d_{R\text{-}T}^{\alpha}(F_X,G_Y)
=
\frac{1}{\alpha-\alpha^{-1}}
\left[
\log I_\alpha
+
I_\alpha
-1
\right],
\end{equation}
where
\begin{equation}
I_\alpha
=
\int_{\lambda}^{\infty}
\int_{x_1}^{\infty}
f_{1,n}^{\alpha}(x_1,x_n)\,
g_{1,n}^{\,1-\alpha}(x_1,x_n)
\,dx_n\,dx_1 ,
\end{equation}
with
$
\lambda=\max(\lambda_1,\lambda_2).
$
Substituting the corresponding joint densities, we obtain
\begin{equation}
\begin{aligned}
I_\alpha
&=
n(n-1)
\beta_1^{2\alpha}
\beta_2^{2(1-\alpha)}
\lambda_1^{2\alpha\beta_1}
\lambda_2^{2(1-\alpha)\beta_2}
\int_{\lambda}^{\infty}
\int_{x_1}^{\infty}
\left[
\left(\frac{\lambda_1}{x_1}\right)^{\beta_1}
-
\left(\frac{\lambda_1}{x_n}\right)^{\beta_1}
\right]^{\alpha(n-2)}
\\
&\quad\times
\left[
\left(\frac{\lambda_2}{x_1}\right)^{\beta_2}
-
\left(\frac{\lambda_2}{x_n}\right)^{\beta_2}
\right]^{(1-\alpha)(n-2)}
(x_1x_n)^{-\left[\alpha(\beta_1+1)+(1-\alpha)(\beta_2+1)\right]}
\,dx_n\,dx_1 .
\end{aligned}
\end{equation}
\end{example}

\begin{example}
\textbf{Rayleigh Distribution}\\

Let
$
X_1,\ldots,X_n
\stackrel{\text{i.i.d.}}{\sim}
\mathrm{Rayleigh}(\beta_1),
\
Y_1,\ldots,Y_n
\stackrel{\text{i.i.d.}}{\sim}
\mathrm{Rayleigh}(\beta_2),
$
where $\beta_1,\beta_2>0$. The corresponding density and distribution functions are
\[
f(x_1)
=
\frac{x_1}{\beta_1^2}
\exp\left(
-\frac{x_1^2}{2\beta_1^2}
\right),
\qquad
F_X(x_1)
=
1-\exp\left(
-\frac{x_1^2}{2\beta_1^2}
\right),
\qquad x_1>0,
\]
and
\[
g(x_1)
=
\frac{x_1}{\beta_2^2}
\exp\left(
-\frac{x_1^2}{2\beta_2^2}
\right),
\qquad
G_Y(x_1)
=
1-\exp\left(
-\frac{x_1^2}{2\beta_2^2}
\right),
\qquad x_1>0.
\]
For the bivariate extreme order statistics $(X_{1:n},X_{n:n})$ and $(Y_{1:n},Y_{n:n})$, the joint densities are
\[
\begin{aligned}
f_{1,n}(x_1,x_n)
&=
n(n-1)
\left[
e^{-x_1^2/(2\beta_1^2)}
-
e^{-x_n^2/(2\beta_1^2)}
\right]^{n-2}
\frac{x_1x_n}{\beta_1^4}
\exp\left(
-\frac{x_1^2+x_n^2}{2\beta_1^2}
\right),
\ 0<x_1<x_n<\infty ,
\end{aligned}
\]
and
\[
\begin{aligned}
g_{1,n}(x_1,x_n)
&=
n(n-1)
\left[
e^{-x_1^2/(2\beta_2^2)}
-
e^{-x_n^2/(2\beta_2^2)}
\right]^{n-2}
\frac{x_1x_n}{\beta_2^4}
\exp\left(
-\frac{x_1^2+x_n^2}{2\beta_2^2}
\right),
\ 0<x_1<x_n<\infty .
\end{aligned}
\]
Then the Rényi--Tsallis divergence of order $\alpha>0$, $\alpha\neq1$ between these two bivariate extreme order statistics is given by
\begin{equation}
d_{R\text{-}T}^{\alpha}(F_X,G_Y)
=
\frac{1}{\alpha-\alpha^{-1}}
\left[
\log I_\alpha
+
I_\alpha
-1
\right],
\end{equation}
where
\begin{equation}
I_\alpha
=
\int_0^\infty
\int_{x_1}^\infty
f_{1,n}^{\alpha}(x_1,x_n)\,
g_{1,n}^{\,1-\alpha}(x_1,x_n)
\,dx_n\,dx_1 .
\end{equation}
Substituting the corresponding joint densities, we obtain
\begin{equation}
\begin{aligned}
I_\alpha
&=
n(n-1)
\beta_1^{-4\alpha}
\beta_2^{-4(1-\alpha)}
\int_0^\infty
\int_{x_1}^\infty
x_1x_n
\left[
e^{-x_1^2/(2\beta_1^2)}
-
e^{-x_n^2/(2\beta_1^2)}
\right]^{\alpha(n-2)}
\\
&\quad\times
\left[
e^{-x_1^2/(2\beta_2^2)}
-
e^{-x_n^2/(2\beta_2^2)}
\right]^{(1-\alpha)(n-2)}
\exp\left[
-\alpha\frac{x_1^2+x_n^2}{2\beta_1^2}
\right]
\\
&\quad\times
\exp\left[
-(1-\alpha)\frac{x_1^2+x_n^2}{2\beta_2^2}
\right]
\,dx_n\,dx_1 .
\end{aligned}
\end{equation}
\end{example}
The above examples have no closed form expression and are evaluated numerically. The following figure illustrates the Rényi--Tsallis divergence for exponential, Weibull, Pareto, and Rayleigh extreme order statistic distributions.
\begin{figure}[h!]
\centering
\begin{minipage}{0.45\textwidth}
\centering
\includegraphics[width=\linewidth,height=4cm,keepaspectratio]{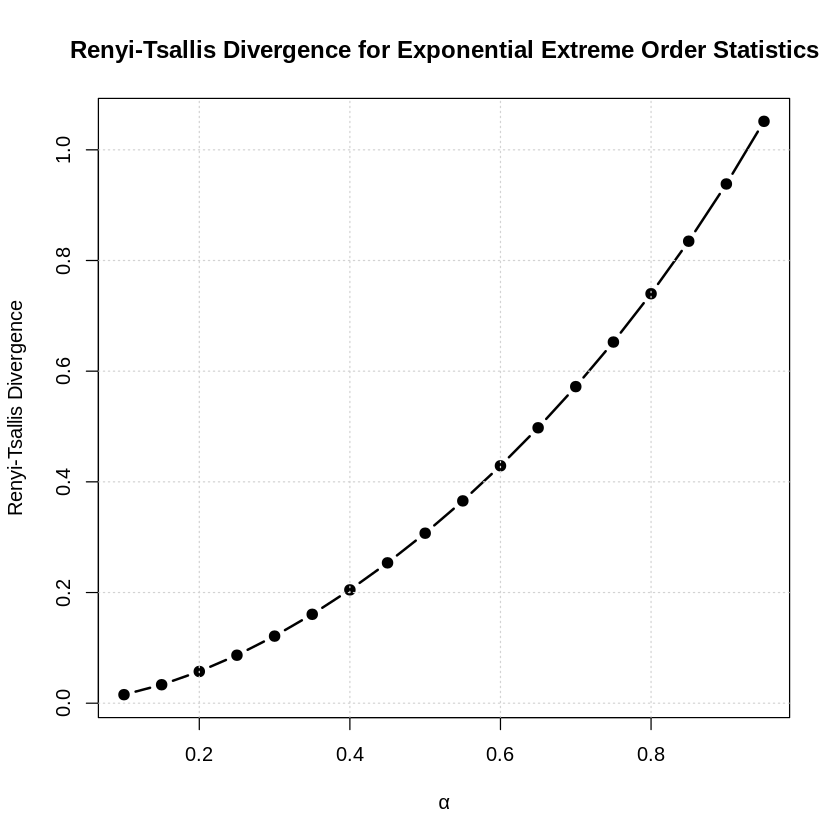}
\end{minipage}
\hfill
\begin{minipage}{0.45\textwidth}
\centering
\includegraphics[width=\linewidth,height=4cm,keepaspectratio]{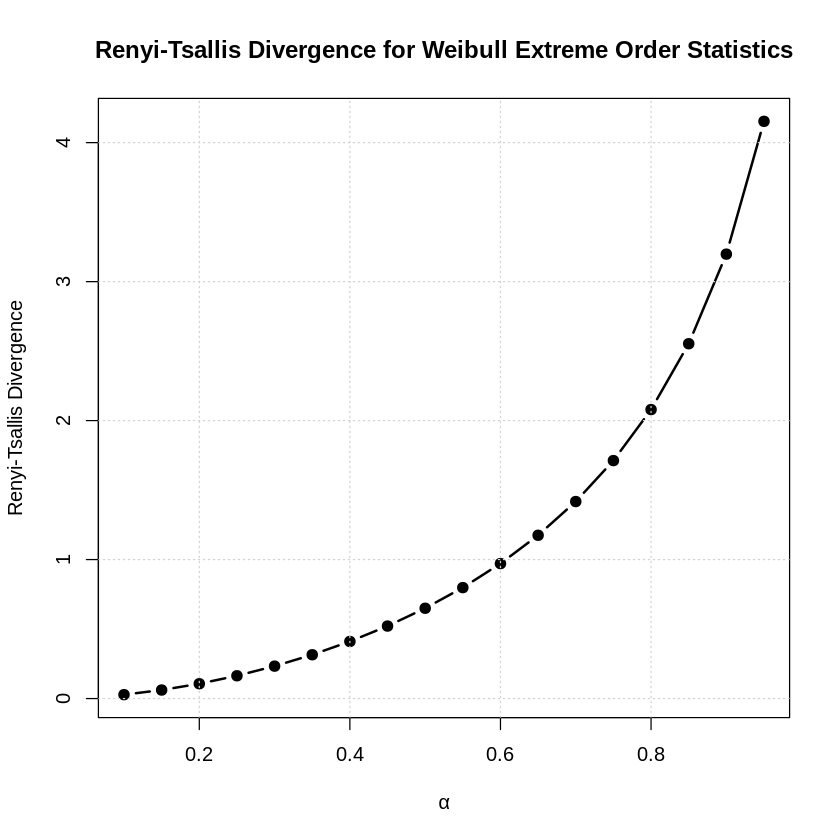}
\end{minipage}

\vspace{0.1cm}

\begin{minipage}{0.45\textwidth}
\centering
\includegraphics[width=\linewidth,height=4cm,keepaspectratio]{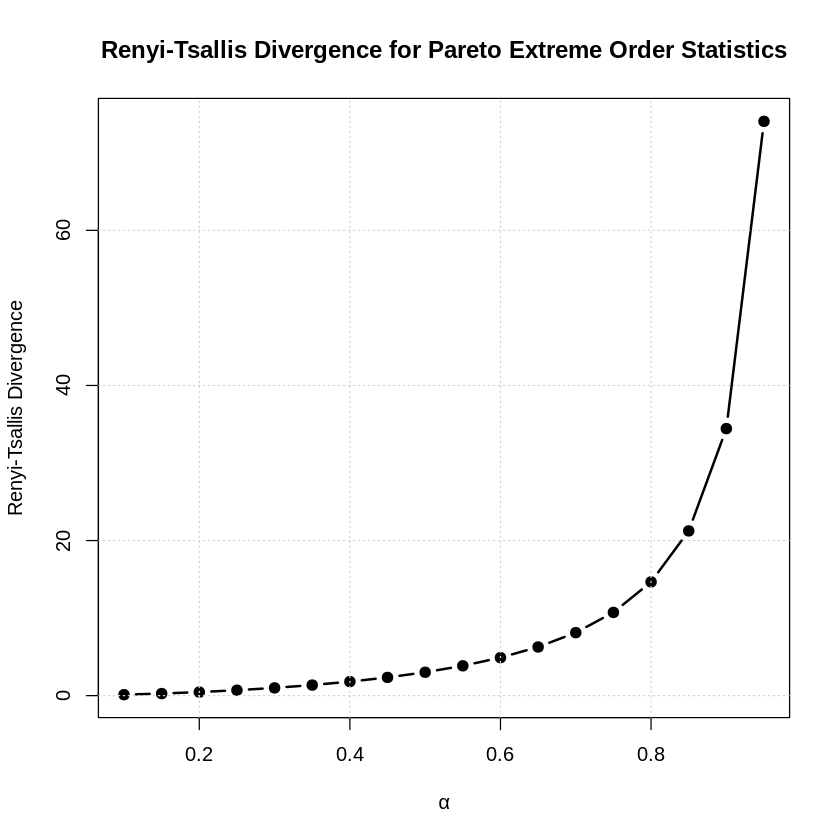}
\end{minipage}
\hfill
\begin{minipage}{0.45\textwidth}
\centering
\includegraphics[width=\linewidth,height=4cm,keepaspectratio]{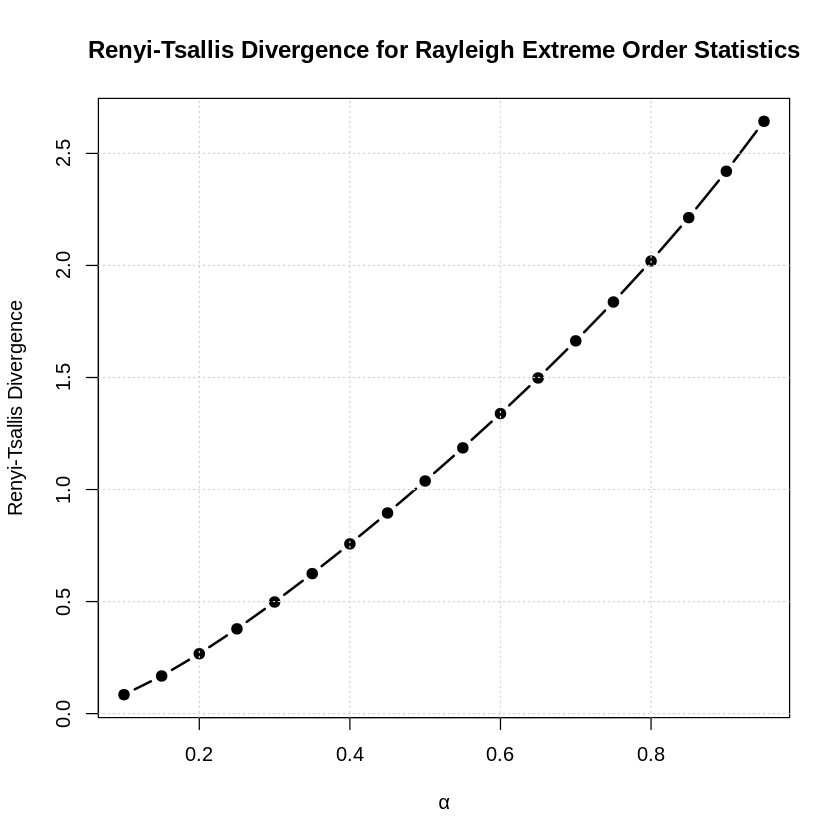}
\end{minipage}

\caption{Rényi--Tsallis divergence for exponential, Weibull, Pareto, and Rayleigh extreme order statistic distributions}
\label{fig:rt_divergence}
\end{figure}
Figure \ref{fig:rt_divergence} shows that the Rényi--Tsallis divergence increases monotonically with $\alpha$ for the extreme order statistics of all four distributions. This indicates that the measure becomes more sensitive to differences between the distributions at larger values of $\alpha$. The increase is most pronounced for the Pareto distribution, particularly for $\alpha > 0.7$, followed by the Weibull distribution. The Rayleigh distribution shows a moderate and steady increase, whereas the exponential distribution exhibits the lowest and smoothest growth.

\section{Nonparametric Estimation of Bivariate Rényi--Tsallis Divergence}
\label{sec:nonparametric}

In this section, nonparametric estimators of the proposed bivariate Rényi--Tsallis divergence measure is developed. Let
$
(X_1,X_2,\ldots,X_{n_1})
\ \text{and} \
(Y_1,Y_2,\ldots,Y_{n_2})
$
be two independent random samples drawn from two populations with unknown bivariate density functions \(f_1(\mathbf{x})\) and \(f_2(\mathbf{x})\), respectively, where
$
\mathbf{x}=(x_1,x_2)^{T}.
$
Here,
$
X_i=(X_{i1},X_{i2})^{T}, \ i=1,2,\ldots,n_1,
$
and
$
Y_j=(Y_{j1},Y_{j2})^{T}, \ j=1,2,\ldots,n_2.$

Since the underlying distributions are unknown, the density functions are estimated nonparametrically using the multivariate kernel density estimator proposed by \cite{duong2007ks}. The kernel density estimators of \(f_1(\mathbf{x})\) and \(f_2(\mathbf{x})\) are given by
\begin{equation}
\begin{aligned}
\widehat{f}_1(\mathbf{x};H_1)
&=
\frac{1}{n_1}
\sum_{i=1}^{n_1}
K_{H_1}(\mathbf{x}-X_i), \\[0.3cm]
\widehat{f}_2(\mathbf{x};H_2)
&=
\frac{1}{n_2}
\sum_{j=1}^{n_2}
K_{H_2}(\mathbf{x}-Y_j).
\end{aligned}
\label{eq:kde_estimators}
\end{equation}
where \(H_1\) and \(H_2\) are symmetric positive definite bandwidth matrices and
\[
K_H(\mathbf{x})
=
|H|^{-1/2}
K\left(H^{-1/2}\mathbf{x}\right).
\]
Throughout this work, the Gaussian kernel
\[
K(\mathbf{x})
=
\frac{1}{(2\pi)^{d/2}}
\exp\left(
-\frac{1}{2}\mathbf{x}^{T}\mathbf{x}
\right)
\]
is employed, where \(d=2\) denotes the dimension of the random vector. The Gaussian kernel is widely used because it provides smooth density estimates and possesses desirable theoretical properties. The bandwidth matrices \(H_1\) and \(H_2\) are selected using the plug-in bandwidth selection method of \cite{duong2007ks}, implemented through the \texttt{Hpi()} function in the \texttt{ks} package in \textsf{R}.

Replacing the unknown densities by their kernel estimates yields the following estimators of the proposed divergence measures.
The estimator of the proposed bivariate Rényi--Tsallis divergence is
\begin{equation}
\widehat{D}_{RT,\alpha}
=
\frac{1}{\alpha-\alpha^{-1}}
\left[
\log
\left(
\int_{\mathbb{R}^{2}}
\widehat{f}_1(\mathbf{x})^{\alpha}
\widehat{f}_2(\mathbf{x})^{1-\alpha}
\, d\mathbf{x}
\right)
+
\left(
\int_{\mathbb{R}^{2}}
\widehat{f}_1(\mathbf{x})^{\alpha}
\widehat{f}_2(\mathbf{x})^{1-\alpha}
\, d\mathbf{x}
-1
\right)
\right].
\label{eq:Rényi_tsallis_estimator}
\end{equation}
where $\alpha>0,\ \alpha\neq1$.
Since these integrals generally do not admit closed-form solutions, they are evaluated numerically over a common two-dimensional grid after normalizing the estimated kernel densities so that they integrate to one.

\subsection*{Application to ACTG315 Longitudinal HIV Data}

The practical applicability of the Rényi--Tsallis divergence measures was demonstrated through the analysis of the ACTG315 longitudinal HIV dataset. The dataset consists of repeated measurements of HIV-1 RNA viral load (RNA) and CD4 lymphocyte counts (CD4) for patients infected with HIV who were receiving antiretroviral therapy. These biomarkers are commonly used to monitor disease progression and treatment response.

The day of measurement was used to classify the observations into three clinical stages: Control, Treatment, and Follow-up. Since multiple observations were available for each patient within a given stage, patient-level averages of RNA and CD4 measurements were calculated. 
For each stage, the joint distribution of RNA viral load and CD4 lymphocyte count was estimated using bivariate Gaussian kernel density estimation with plug-in bandwidth selection.

Pairwise Rényi, Tsallis, and Rényi-Tsallis divergence measures were calculated between the Control–Treatment, Control–Follow-up, and Treatment–Follow-up comparisons, using the estimated densities for $\alpha=0.4$. The corresponding values are presented in Tables~\ref{tab:Rényi_divergence},
\ref{tab:tsallis_divergence}, and \ref{tab:Rényi_tsallis_divergence}. Results showed clear distributional differences between the three clinical stages. The biggest divergence values were observed between the Control and Follow-up groups. The Treatment and Follow-up groups had relatively smaller divergence values, indicating more similarity in their joint distributions of biomarkers. Also, out of the three divergence measures Rényi--Tsallis divergence gives the smallest values for the pairwise comparison.

To assess uncertainty of the estimated divergences, patient-level nonparametric bootstrap (1000 replications) was performed. Kernel densities and divergence measures were recomputed for each bootstrap sample and bootstrap bias, variance, standard error, mean squared error (MSE) and 95\% bootstrap percentile confidence intervals were obtained.
 
Table \ref{tab:estimator_performance} presents the bootstrap comparison of the estimators of Rényi, Tsallis and Rényi-Tsallis divergence for bias, variance, standard error, mean squared error (MSE) and confidence interval. For both control versus treatment and treatment versus follow-up Rényi--Tsallis divergence has the lowest variance, standard error and MSE. So for most of the pairwise clinical stage comparison we can say that the Rényi--Tsallis divergence estimator showed lower variability and MSE, suggesting greater numerical stability.

\begin{table*}[t]
\centering
\small
\setlength{\tabcolsep}{3pt}

\begin{minipage}[t]{0.48\textwidth}
\centering
\tbl{Rényi divergence values between the Control, Treatment, and Follow-up stages}
{\begin{tabular}{lccc}
\toprule
& Control & Treatment & Follow-up \\
\midrule
Control   & 0.0000 & 1.0795 & 1.5684 \\
Treatment & 1.1214 & 0.0000 & 0.3819 \\
Follow-up & 1.2540 & 0.3433 & 0.0000 \\
\bottomrule
\end{tabular}}
\label{tab:Rényi_divergence}
\end{minipage}
\hfill
\begin{minipage}[t]{0.48\textwidth}
\centering
\tbl{Tsallis divergence values between the Control, Treatment, and Follow-up stages}
{\begin{tabular}{lccc}
\toprule
& Control & Treatment & Follow-up \\
\midrule
Control   & 0.0000 & 0.7946 & 1.0163 \\
Treatment & 0.8163 & 0.0000 & 0.3413 \\
Follow-up & 0.8813 & 0.3102 & 0.0000 \\
\bottomrule
\end{tabular}}
\label{tab:tsallis_divergence}
\end{minipage}

\vspace{0.5cm}

\begin{minipage}[t]{0.48\textwidth}
\centering
\tbl{Rényi--Tsallis divergence values between the Control, Treatment, and Follow-up stages}
{\begin{tabular}{lccc}
\toprule
& Control & Treatment & Follow-up \\
\midrule
Control   & 0.0000 & 0.5355 & 0.7385 \\
Treatment & 0.5536 & 0.0000 & 0.2066 \\
Follow-up & 0.6101 & 0.1867 & 0.0000 \\
\bottomrule
\end{tabular}}
\label{tab:Rényi_tsallis_divergence}
\end{minipage}

\end{table*}

\begin{table}
\tbl{Bootstrap bias, variance, standard error, MSE, and 95\% percentile confidence intervals of the divergence estimators for pairwise clinical-stage comparisons}
{\begin{tabular}{llccccc}
\toprule
\textbf{Comparison} & \textbf{Estimator} & \textbf{Bias} & \textbf{Variance} & \textbf{SE} & \textbf{MSE} & \textbf{95\% CI} \\
\midrule
Control vs Treatment
& Rényi & 0.0249 & 0.03688 & 0.1920 & 0.03750 & (0.8261, 1.5667) \\
& Tsallis & 0.0075 & 0.00872 & 0.0934 & 0.00878 & (0.6514, 1.0156) \\
& Rényi--Tsallis & 0.0092 & 0.00664 & 0.0815 & 0.00672 & (0.4221, 0.7378) \\
\midrule
Control vs Follow-up
& Rényi & 0.0172 & 0.07206 & 0.2684 & 0.07236 & (1.1967, 2.2442) \\
& Tsallis & -0.0012 & 0.00914 & 0.0956 & 0.00914 & (0.8538, 1.2331) \\
& Rényi--Tsallis & 0.0046 & 0.01078 & 0.1038 & 0.01080 & (0.5858, 0.9935) \\
\midrule
Treatment vs Follow-up
& Rényi & 0.0380 & 0.01404 & 0.1185 & 0.01548 & (0.2411, 0.6861) \\
& Tsallis & 0.0267 & 0.00810 & 0.0900 & 0.00881 & (0.2245, 0.5624) \\
& Rényi--Tsallis & 0.0185 & 0.00355 & 0.0595 & 0.00389 & (0.1330, 0.3567) \\
\bottomrule
\end{tabular}}
\label{tab:estimator_performance}
\end{table}

\section{Conclusion}
\label{sec:conclusion}

In this article, we have proposed bivariate Rényi--Tsallis divergence, a natural extension for the univariate Rényi--Tsallis divergence which combines features of Rényi and Tsallis measures in a single framework. Several properties of bivariate Rényi--Tsallis divergence including non--negativity, pseudo additivity, behavior under linear and monotone transformations are studied. A nonparametric kernel based estimator was also developed in cases where the underlying distributions are unknown. The methodology has been successfully applied to ACTG315 longitudinal HIV data set. Bootstrap analysis has showed that the proposed estimator is able to provide accurate and stable divergence estimates. The present work can be extended to divergence measures based on bivariate distribution functions, bivariate survival functions and bivariate cumulative residual distributions. In addition, the proposed Rényi--Tsallis divergence may be generalized to the multivariate setting for comparing higher dimensional random vectors.

\section*{Declarations}
\textbf{Conflict of interest} On behalf of all authors, the corresponding author states that there is no conflict of
interest.

\bibliographystyle{plainnat}
\bibliography{references}

@ARTICLE{rafflesia2026,
author = {Rafflesia Chackochan, M. and Sankaran, P. G. and Unnikrishnan Nair, N.},
year = 2026,
title = {Bivariate {K}ullback--{L}eibler Divergence},
journal = {Communications in Statistics -- Theory and Methods},
volume = 55,
number = 1,
pages = {292--312}
}

@INCOLLECTION{clausius1865,
author = {Clausius, R.},
year = 1865,
title = {On Different Forms of the Fundamental Equations of the Mechanical Theory of Heat and Their Convenience for Application},
booktitle = {The Second Law of Thermodynamics},
pages = {162--193}
}

@ARTICLE{khammar2022,
author = {Khammar, A. and Noughabi, H. A.},
year = 2022,
title = {On Differential {R}enyi's--{T}sallis Divergence Measure and Its Applications},
journal = {Journal of Statistical Computation and Simulation},
volume = 92,
number = 15,
pages = {3109--3128},
doi = {10.1080/00949655.2022.2043417}
}

@ARTICLE{kullback1951,
author = {Kullback, S. and Leibler, R. A.},
year = 1951,
title = {On Information and Sufficiency},
journal = {The Annals of Mathematical Statistics},
volume = 22,
number = 1,
pages = {79--86},
doi = {10.1214/aoms/1177729694}
}

@INPROCEEDINGS{renyi1961,
author = {Renyi, A.},
year = 1961,
title = {On Measures of Entropy and Information},
editor = {Neyman, J.},
booktitle = {Proceedings of the Fourth Berkeley Symposium on Mathematical Statistics and Probability, Volume 1: Contributions to the Theory of Statistics},
address = {Berkeley (CA)},
publisher = {University of California Press},
pages = {547--562}
}

@ARTICLE{shannon1948,
author = {Shannon, C. E.},
year = 1948,
title = {A Mathematical Theory of Communication},
journal = {The Bell System Technical Journal},
volume = 27,
number = 3,
pages = {379--423},
doi = {10.1002/j.1538-7305.1948.tb01338.x}
}

@ARTICLE{tsallis1988,
author = {Tsallis, C.},
year = 1988,
title = {Possible Generalization of {B}oltzmann--{G}ibbs Statistics},
journal = {Journal of Statistical Physics},
volume = 52,
number = {1--2},
pages = {479--487},
doi = {10.1007/BF01016429}
}

@ARTICLE{wondie2017joint,
author = {Wondie, L. and Kumar, S.},
year = 2017,
title = {A Joint Representation of {R}enyi's and {T}sallis's Entropy with Application in Coding Theory},
journal = {International Journal of Mathematics and Mathematical Sciences},
volume = 2017,
eid = {2683293},
doi = {10.1155/2017/2683293}
}

@ARTICLE{perlman1974jensen,
author = {Perlman, Michael D.},
year = 1974,
title = {Jensen's Inequality for a Convex Vector-Valued Function on an Infinite-Dimensional Space},
journal = {Journal of Multivariate Analysis},
volume = 4,
number = 1,
pages = {52--65},
doi = {10.1016/0047-259X(74)90020-8}
}

@ARTICLE{dragomir2000,
author = {Dragomir, S. S. and Gluscevic, V.},
year = 2000,
title = {Some Inequalities for the {K}ullback--{L}eibler and {$\chi^2$}-Distances in Information Theory and Applications},
journal = {RGMIA Research Report Collection},
volume = 3,
number = 2,
pages = {199--210}
}

@ARTICLE{bental1991,
author = {Ben-Tal, Aharon and Ben-Israel, Adi and Teboulle, Marc},
year = 1991,
title = {Certainty Equivalents and Information Measures: Duality and Extremal Principles},
journal = {Journal of Mathematical Analysis and Applications},
volume = 157,
number = 1,
pages = {211--236},
doi = {10.1016/0022-247X(91)90140-N}
}

@BOOK{cichocki2009,
author = {Cichocki, Andrzej and Zdunek, Rafal and Phan, Anh Huy and Amari, Shun-ichi},
year = 2009,
title = {Nonnegative Matrix and Tensor Factorizations: Applications to Exploratory Multi-Way Data Analysis and Blind Source Separation},
address = {Chichester (UK)},
publisher = {John Wiley \& Sons}
}

@ARTICLE{bhattacharyya1943,
author = {Bhattacharyya, Anil},
year = 1943,
title = {On a Measure of Divergence Between Two Statistical Populations Defined by Their Probability Distributions},
journal = {Bulletin of the Calcutta Mathematical Society},
volume = 35,
pages = {99--110}
}

@BOOK{hellinger1907,
author = {Hellinger, Ernst},
year = 1907,
title = {Die Orthogonalinvarianten Quadratischer Formen von Unendlichvielen Variabelen},
address = {Leipzig (DE)},
publisher = {W. Fr. Kaestner}
}

@MISC{duong2007ks,
author = {Duong, Tarn},
year = 2007,
title = {Using \texttt{ks} for Bivariate Kernel Density Estimation},
note = {R package documentation}
}

@ARTICLE{Ebrahimi2007,
  author  = {Ebrahimi, Nader and Kirmani, S. N. U. A. and Soofi, Ehsan S.},
  title   = {Multivariate Dynamic Information},
  journal = {Journal of Multivariate Analysis},
  volume  = {98},
  number  = {2},
  pages   = {328--349},
  year    = {2007},
  doi     = {10.1016/j.jmva.2005.08.004}
}

@inproceedings{zhang2023properties,
  author    = {Zhang, Yufeng and Pan, Jialu and Li, Kenli and Liu, Wanwei and
               Chen, Zhenbang and Liu, Xinwang and Wang, Ji},
  title     = {On the Properties of {K}ullback--{L}eibler Divergence Between
               Multivariate {G}aussian Distributions},
  booktitle = {Advances in Neural Information Processing Systems},
  volume    = {36},
  pages     = {58152--58165},
  year      = {2023}
}

@article{contrerasreyes2012,
  author  = {Contreras-Reyes, Javier E. and Arellano-Valle, Reinaldo B.},
  title   = {Kullback--Leibler Divergence Measure for Multivariate Skew-Normal Distributions},
  journal = {Entropy},
  volume  = {14},
  number  = {9},
  pages   = {1606--1626},
  year    = {2012},
  doi     = {10.3390/e14091606}
}

\end{document}